\documentclass[11pt]{amsart}

\usepackage{amsmath,amsthm,amssymb}
\usepackage{hhline}
\usepackage{longtable,array}
\usepackage{mathrsfs}

\newtheorem{Le}{Lemma}[section]

\newtheorem{Th}[Le]{Theorem}
\newtheorem{Co}[Le]{Corollary}

\theoremstyle{definition}
\newtheorem{Def}[Le]{Definition}

\newtheorem{Conv}[Le]{Convention}

\theoremstyle{remark}
\newtheorem{Ex}[Le]{Example}
\newtheorem{Rem}[Le]{Remark}

\allowdisplaybreaks
\newcommand{\Kthree}{\mathscr{K}\!\mathit{3}}

\begin{document}

\title[K3 surfaces with a pair of involutions]{K3 surfaces with a pair of 
commuting non-symplectic involutions}
\author{Frank Reidegeld}
\subjclass[2010]{14J28, 53C26} 

\begin{abstract}
We study K3 surfaces with a pair of commuting involutions that are non-symplectic
with respect to two anti-commuting complex structures that are determined by a 
hyper-K\"ahler metric. One motivation for this paper is the role of such 
$\mathbb{Z}^2_2$-actions for the construction of $G_2$-manifolds. We find a 
large class of smooth K3 surfaces with such pairs of involutions, but we also 
pay special attention to the case that the K3 surface has ADE-singularities. Therefore, 
we introduce a special class of non-symplectic involutions that are suitable for explicit 
calculations and find $320$ examples of pairs of involutions that act on K3 surfaces 
with a great variety of singularities.  
\end{abstract}

\maketitle

\setcounter{tocdepth}{1}
\tableofcontents

\section{Introduction} 

A non-symplectic involution of a K3 surface $S$ is a holomorphic involution 
$\rho:S\rightarrow S$ such that $\rho$ acts as $-1$ on $H^{2,0}(S)$. Any
K3 surface with a non-symplectic involution admits a K\"ahler metric that 
is invariant under $\rho$. Since any K\"ahler metric on a K3 surface is in
fact hyper-K\"ahler, there are three complex structures $I$, $J$ and $K$ 
and three K\"ahler forms $\omega_I$, $\omega_J$ and $\omega_K$ on 
$S$. If $\rho$ is holomorphic with respect to $I$, we have

\[
\rho^{\ast}\omega_I = \omega_I\:, 
\qquad 
\rho^{\ast}\omega_J = - \omega_J\:,
\qquad
\rho^{\ast}\omega_K = -\omega_K 
\]   

In this paper, we search for K3 surfaces that admit two non-symplectic involutions
$\rho^1$ and $\rho^2$. We require that $\rho^1$ and $\rho^2$ commute and that
they are non-symplectic with respect to two different complex structures from the
triple $(I,J,K)$. Without loss of generality, we can assume that 

\begin{equation}
\label{KahlerRelations1}
\begin{aligned}
& {\rho^1}^{\ast}\omega_I = \omega_I\:, 
\qquad 
{\rho^1}^{\ast}\omega_J = - \omega_J\:,
\qquad
{\rho^1}^{\ast}\omega_K = -\omega_K \\
& {\rho^2}^{\ast}\omega_I = - \omega_I\:, 
\qquad 
{\rho^2}^{\ast}\omega_J = \omega_J\:,
\qquad
{\rho^2}^{\ast}\omega_K = - \omega_K 
\end{aligned}
\end{equation}

A motivation to study these pairs of involutions is their relation to the 
construction of $G_2$-manifolds. $\rho^1$ and $\rho^2$ generate a group 
that is isomorphic to $\mathbb{Z}^2_2$ and acts isometrically on $S$. 
In the habilitation thesis of the author \cite{ReidHabil} we have described 
how such an action can be extended to products $S\times T^3$ of a K3 
surface and a $3$-torus such that the quotients $(S\times T^3)/\mathbb{Z}^2_2$
carry a $G_2$-structure. In a forthcoming paper we resolve the singularities
of those quotients by the methods of Karigiannis and Joyce \cite{KarJoy} and 
thus obtain compact $G_2$-manifolds. In \cite{ReidHabil} we have also 
explained that pairs of involutions with the above properties can be used to solve 
the so called matching problem in Kovalev's and Lee's construction of compact 
$G_2$-manifolds by twisted connected sums \cite{KoLe}.

A non-symplectic involution is determined by its action on the lattice \linebreak
$H^2(S,\mathbb{Z})$. Nikulin \cite{Nikulin1,Nikulin2,Nikulin3} has classified 
the non-symplectic involutions of K3 surfaces in terms of invariants of their
fixed lattices. By embedding the direct sum of two possible fixed lattices into 
$H^2(S,\mathbb{Z})$ we are able to find a large class of pairs $(\rho^1,\rho^2)$ 
with the desired properties. If the hyper-K\"ahler metric on $S$ is chosen
generically, it has no singularities. In this paper, we also pay special attention to 
the case that $S$ has ADE-singularities. If we start one of the above constructions of 
$G_2$-manifolds with a K3 surface with singularities as its starting point, 
we obtain a $G_2$-orbifold with ADE-singularities. Such orbifolds are studied as
compactifications of M-theory since the ADE-singularities are needed to explain
the presence of non-abelian gauge fields \cite{Ach,AchGuk}.   

In order to construct K3 surfaces with a pair $(\rho^1,\rho^2)$ satisfying 
(\ref{KahlerRelations1}) that have as many types of ADE-singularities
as possible, we restrict ourselves to a special class of non-symplectic
involutions that are suitable for explicit calculations. We call these
involutions simple. We show that $28$ out of $75$ types of
non-symplectic involutions are simple. Moreover, we find $320$ 
different kinds of pairs $(\rho^1,\rho^2)$ such that $\rho^1$ and 
$\rho^2$ are simple. Each of them acts on a K3 surface with 
$3$ $A_1$- and $2$ $E_8$-singularities. Furthermore, there exist 
plenty of K3 surfaces with fewer and milder singularities 
that admit the same kind of involutions. 

This paper is organised as follows. In Section \ref{K3Intro} - \ref{NonsympSection}, 
we present the necessary background material that can be found in the literature. 
The definition of a simple non-symplectic involution and their classification can be 
found at the end of Section \ref{NonsympSection}. The main results of this paper are 
proven in Section \ref{OneInvolution} and \ref{TwoInvolutions}. Section \ref{OneInvolution} 
deals with K3 surfaces with singularities that admit one involution and Section 
\ref{TwoInvolutions} deals with K3 surfaces that admit a pair of involutions.

\section{K3 surfaces and their moduli spaces}
\label{K3Intro}

Since some of the readers of this article may have a background in
Riemannian rather than algebraic geometry, we provide a short
introduction to the theory of K3 surfaces and their moduli spaces. We 
refer the reader to \cite[Chapter VIII]{BHPV} and references therein
for a more detailed account. First of all, we define what a K3 surface is. 

\begin{Def}
A \emph{K3 surface} is a compact, simply connected, complex surface with 
trivial canonical bundle.  
\end{Def}

At the beginning of this section, we consider only smooth K3 surfaces. Later on,
we allow ADE-singularities, too. The underlying manifold of any K3 surface is of
a fixed diffeomorphism type. Therefore, the topological invariants of all K3 
surfaces are the same. The second  cohomology with integer coefficients together
with the intersection form is a lattice. Since  we have to work with the second 
cohomology and its various sublattices, we need some concepts from lattice theory. The 
content of the following pages can be found in any reference on this subject, for
example in \cite[Chapter I.2]{BHPV}, \cite{Do} or in \cite{Ebeling}. 

\begin{Def}
\begin{enumerate}
    \item A \emph{lattice} is a free abelian group $L$ of finite rank together with a
    symmetric bilinear form $\cdot:L\times L\rightarrow \mathbb{Z}$. We write $x^2$ 
    for $x\cdot x$. The \emph{rank} of a lattice is the same as the rank of the underlying
    abelian group. $L$ is called \emph{even} if $x^2 \in 2\mathbb{Z}$ for all $x\in L$.
    Let $(e_1,\ldots,e_n)$ be a basis of $L$. The $n\times n$-matrix with coefficients 
    $e_i\cdot e_j$ is called the \emph{Gram matrix $G(L)$ of $L$ with respect to 
   $(e_1,\ldots,e_n)$}. $L$ is called \emph{unimodular} if $|\det{G(L)}| = 1$. 
    
    \item The $\mathbb{Z}$-bilinear form on $L$ can be extended to an $\mathbb{R}$-bilinear 
    form on $L\otimes_{\mathbb{Z}} \mathbb{R}$. Terms as \emph{non-degenerate 
    lattice} and \emph{signature of a lattice} will always be defined with respect to the 
    extended form.     
         
    \item Let $L$ and $L'$ be lattices and let $\cdot_L$ and $\cdot_{L'}$ be the 
    corresponding bilinear forms. A \emph{lattice isomorphism of $L$ and $L'$}\: 
    is a bijective $\mathbb{Z}$-linear map $\phi:L\rightarrow L'$ with $x\cdot_L y = 
    \phi(x)\cdot_{L'} \phi(y)$ for all $x,y\in L$. If $L=L'$, $\phi$ is called an
    \emph{automorphism}. We denote the group of all automorphisms of $L$
    by $\text{Aut}(L)$. 
\end{enumerate}         
\end{Def}

\begin{Rem}
The number $|\det{G(L)}|$ from the above definition
is independent of the choice of the basis $(e_1,\ldots,e_n)$. 
\end{Rem}

\begin{Def}
\begin{enumerate}
    \item An element $x$ of a lattice $L$ is called \emph{primitive} if there exists 
    no $k>1$ and $y\in L$ such that $x=k\cdot y$.

    \item A sublattice $K\subseteq L$ is called \emph{primitive} if the quotient
    $L/K$ has no torsion. 
    
    \item A lattice $N$ is \emph{primitively embedded} in $L$ if $L$ has a primitive 
    sublattice that is isomorphic to $N$.       
\end{enumerate} 
\end{Def}   

The dual of a lattice $L$ is defined as

\begin{equation*}
L^{\ast} := \{\phi:L\rightarrow \mathbb{Z} | \phi\:\:\text{is $\mathbb{Z}$-linear} \}\:.
\end{equation*} 

From now on, we assume that $L$ is a non-degenerate lattice. $L^{\ast}$ can 
be equipped with the dual bilinear form, which takes its values in $\mathbb{Q}$
but not necessarily in $\mathbb{Z}$. The Gram matrix of $L^{\ast}$ with respect
to the dual basis is given by $G(L)^{-1}$. If $L$ is unimodular, $L^{\ast}$ is thus 
a lattice, too. The map $\imath: L\rightarrow L^{\ast}$ that is defined by 
$\imath(x)(y) := x\cdot y$ is an injection. The quotient group $L^{\ast}/\imath(L)$
is called the \emph{discriminant group of $L$}. 

\begin{Le}
The discriminant group of a lattice $L$ is a finite group of order $|\det{G(L)}|$. 
The minimal number $\ell(L)$ of generators of the discriminant group satisfies
$\ell(L)\leq \text{rank}\:L$.  
\end{Le}

The invariant $\ell(L)$ allows us to formulate a theorem on primitive embeddings 
that can be found in \cite{Do} or \cite{Nikulin2}.

\begin{Th}
\label{LatticeEmbThm}
Let $K$ be an even non-degenerate lattice of signature $(k_+,k_-)$ and
$L$ be an even unimodular lattice of signature $(l_+,l_-)$. We assume that
$k_+\leq l_+$ and $k_-\leq l_-$ and that 

\begin{enumerate}
    \item $2\cdot\text{rank}(K) \leq \text{rank}(L)$ or
    
    \item $\text{rank}(K) + \ell(K) < \text{rank}(L)$.
\end{enumerate}

Then there exists a primitive embedding $i:K\rightarrow L$. If in addition 
$k_+ < l_+$ and $k_- < l_-$ and one of the following conditions holds

\begin{enumerate}
    \item  $2\cdot\text{rank}(K) \leq \text{rank}(L) - 2$, 
    
    \item $\text{rank}(K) + \ell(K) \leq \text{rank}(L) - 2$,
\end{enumerate}

the embedding $i$ is unique up to an automorphism of $L$. 
\end{Th}

We return to K3 surfaces and describe their topology.

\begin{Th}
Let $S$ be a K3 surface.
\begin{enumerate}
    \item The Hodge numbers of $S$ are determined by 
    $h^{0,0}(S)= h^{2,0}(S)=1$, $h^{1,0}(S)=0$ and $h^{1,1}(S)=20$.  

    \item The second integral cohomology $H^2(S,\mathbb{Z})$ 
    together with the intersection form is an even unimodular lattice
    of signature $(3,19)$. Up to isometries, the only lattice with these
    properties is  

    \begin{equation*}
    L := 3 H \oplus 2(-E_8) \:,
    \end{equation*}

    where $H$ is the hyperbolic plane lattice with the bilinear form

    \begin{equation}
    \label{Hyperbolic}
    \left(\,
    \begin{array}{cc}
    0 & 1 \\
    1 & 0 \\
    \end{array}
    \,\right)
    \end{equation}

   and $-E_8$ is the root lattice of $E_8$ together with the negative of the usual 
   bilinear form.
\end{enumerate}  
\end{Th}

These facts motivate the following definitions. 

\begin{Def}
\begin{enumerate}
    \item The lattice $L$ from the above theorem is called the \emph{K3 lattice}.
    
    \item A K3 surface $S$ together with a lattice isometry $\phi:H^2(S,\mathbb{Z})
    \rightarrow L$ is called a \emph{marked K3 surface}.
   
    \item Two marked K3 surfaces $(S,\phi)$ and $(S',\phi')$ are called
    \emph{isomorphic} if there exists a biholomorphic map $f:S\rightarrow S'$
    such that $\phi\circ f^{\ast} = \phi'$, where $f^{\ast}: H^2(S',\mathbb{Z})
    \rightarrow H^2(S,\mathbb{Z})$ is the pull-back. 
\end{enumerate}
\end{Def}

The first Chern class on $S$ is a bijective map between the Picard group and 
$H^{1,1}(S)\cap H^2(S,\mathbb{Z})$. Therefore, we introduce the following
terms:

\begin{Def}
The lattice $H^{1,1}(S)\cap H^2(S,\mathbb{Z})$ is called the 
\emph{Picard lattice} and its rank is called the \emph{Picard number}. The 
orthogonal complement of the Picard lattice in $H^2(S,\mathbb{Z})$ is 
called the \emph{transcendental lattice}.  
\end{Def}

\begin{Conv}
The maximal value of the Picard number is $20$. In the literature, a K3
surface with maximal Picard number is often called singular and a compact,
simply connected, complex surface with trivial canonical bundle that may 
admit ADE-singularities is sometimes called a \emph{Gorenstein K3 surface}.
In this article, we use a different convention and call K3 surfaces with 
ADE-singularities singular.   
\end{Conv}

Any K3 surface $S$ admits a K\"ahler metric. Since $S$ has trivial canonical 
bundle, there exists a unique Ricci-flat K\"ahler metric in each K\"ahler class. 
The holonomy group $SU(2)$ is isomorphic to $Sp(1)$. Therefore, the 
Ricci-flat K\"ahler metrics are in fact hyper-K\"ahler. When we talk about 
isomorphisms between K3 surfaces, we usually mean biholomorphic maps
with respect to fixed complex structures on the K3 surfaces. 
Another natural class of maps between K3 surfaces, which will be studied later
on, are isometries between K3 surfaces with hyper-K\"ahler metrics. It should 
be noted that there are isometries between K3 surfaces that are not holomorphic. 

We need some background knowledge on the moduli spaces of K3 surfaces. There
are several related moduli spaces whose points represent K3 surfaces with an extra 
structure. We denote them all by $\Kthree$ with an appropriate index. The first of them
is the \emph{moduli space of marked K3 surfaces} $\Kthree^m$ that is defined as the
set of all marked K3 surfaces modulo isomorphisms. We describe $\Kthree^m$ in more 
detail below. 

On any K3 surface, there exists a holomorphic $(2,0)$-form that
is unique up to multiplication with a constant. We denote it by 
$\omega_J + i\omega_K$, where $\omega_J$ and $\omega_K$ are 
real-valued 2-forms. This observation motivates the following definition. 

\begin{Def}
Let $(S,\phi)$ be a marked K3 surface. Moreover, let $\mathbb{K}\in
\{\mathbb{R},\mathbb{C}\}$, $L_{\mathbb{K}}:= L \otimes_{\mathbb{Z}} 
\mathbb{K}$ and $\phi_{\mathbb{K}}: H^2(S,\mathbb{K}) 
\rightarrow L_{\mathbb{K}}$ be the $\mathbb{K}$-linear 
extension of $\phi$. The complex line that is spanned by 
$\phi_{\mathbb{C}} ([\omega_J + i\omega_K])$, where the square
brackets denote the cohomology class, defines a point 
$p(S,\phi)\in \mathbb{P}(L_{\mathbb{C}})$, where 
$\mathbb{P}(L_{\mathbb{C}})$ is the projective space of all
complex lines in $L_{\mathbb{C}}$. $p(S,\phi)$ is called the 
\emph{period point of $(S,\phi)$}. This assignment defines a map 
$p:\Kthree^m\rightarrow \mathbb{P}(L_{\mathbb{C}})$, which is called 
the \emph{period map for K3 surfaces}. 
\end{Def}

It is not difficult to prove that $p(S,\phi)$ is always contained in the following
subset of $\mathbb{P}(L_{\mathbb{C}})$. 

\begin{Def}
We denote the complex line that is spanned by $x\in
L_{\mathbb{C}}\setminus\{0\}$ by $\ell_x$. 
The set 
    
\begin{equation*}
\Omega:=\{ \ell_x \in \mathbb{P}(L_{\mathbb{C}}) | x\cdot x = 0,\: x\cdot\overline{x}>0 \}  
\end{equation*}
    
is called the \emph{period domain}. 
\end{Def} 

We reduce the target set of the period map such that from now on  
$p:\Kthree^m\rightarrow \Omega$. An important theorem
in the theory of K3 surfaces is that the period map is surjective. 
Moreover, it is a local isomorphism of complex manifolds, but
it is not injective. Therefore, $\Kthree^m$ and $\Omega$ are
not isomorphic. In order to describe $\Kthree^m$ explicitly, we need some 
further definitions. 

\begin{Def}
\begin{enumerate}
    \item Let $S$ and $S'$ be K3 surfaces. A lattice isometry 
    $\psi:H^2(S,\mathbb{Z}) \rightarrow H^2(S',\mathbb{Z})$ is called a 
    \emph{Hodge-isometry} if its $\mathbb{C}$-linear extension preserves 
    the Hodge decomposition $H^2(S,\mathbb{C})=H^{2,0}(S) \oplus H^{1,1}(S) 
    \oplus H^{0,2}(S)$. 

    \item A class $x\in H^2(S,\mathbb{Z})$ is called \emph{effective} if
    there exists an effective divisor $D$ of $S$ with $c_1(\mathcal{O}_S(D))
    = x$. An effective class $x$ is called \emph{nodal} if $x^2 = -2$. 

   \item The connected component of the set $\{x\in H^{1,1}(S,\mathbb{R}) | x\cdot x > 
   0\}$ which contains a K\"ahler class is called the \emph{positive
   cone} of $S$.

    \item A Hodge-isometry $\psi:H^2(S,\mathbb{Z}) \rightarrow 
    H^2(S',\mathbb{Z})$ is called \emph{effective} if it maps the 
    positive cone of $S$ to the positive cone of $S'$ and effective 
    classes in $H^2(S,\mathbb{Z})$ to effective classes in 
    $H^2(S',\mathbb{Z})$. 
\end{enumerate}
\end{Def}

\begin{Rem}
Since $H^{1,1}(S) = \overline{H^{1,1}(S)}$, $H^{1,1}(S)$ is a complex
vector space. We denote its real part $H^{1,1}(S)\cap H^2(S,\mathbb{R})$
by $H^{1,1}(S,\mathbb{R})$. The restriction of the intersection form to 
$H^{1,1}(S,\mathbb{R})$ has signature $(1,19)$. The set $\{x\in 
H^{1,1}(S,\mathbb{R}) | x\cdot x > 0\}$ thus has exactly two connected 
components. Exactly one of them contains a K\"ahler class and the 
definition of the positive cone therefore makes sense.
\end{Rem} 

The following lemma often helps to decide if a Hodge-isometry is effective.

\begin{Le} (See \cite[p. 313]{BHPV})
\label{EffectiveLemma}
Let $S$ and $S'$ be K3 surfaces and $\psi:H^2(S,\mathbb{Z})
\rightarrow H^2(S',\mathbb{Z})$ be a Hodge-isometry.  If $\psi$
maps at least one K\"ahler class of $S$ to a K\"ahler class of
$S'$, then $\psi$ is effective.
\end{Le}

With help of the terms that we have defined above, we are able to
state the following theorem. 

\begin{Th} (Torelli theorem)
Let $S$ and $S'$ be two unmarked K3 surfaces. If there exists an effective
Hodge-isometry $\psi:H^2(S',\mathbb{Z}) \rightarrow H^2(S,\mathbb{Z})$,
$\psi$ is the pull-back of a unique biholomorphic map $f:S\rightarrow S'$. 
\end{Th}

The converse of the above theorem is also true. If $f:S\rightarrow S'$ 
is a biholomorphic map, its pull-back is an effective Hodge-isometry.
Since effective Hodge-isometries are closely related to the set of 
all K\"ahler classes,  we are able to describe $\Kthree^m$ with help of
an explicit description of the K\"ahler cone. 

\begin{Th}
Let $S$ be a K3 surface and let $\mathcal{C}_S\subseteq H^{1,1}(S,\mathbb{R})$ be 
its K\"ahler cone, i.e. the set of all cohomology classes representing
a K\"ahler form. Then we have

\begin{equation*}
\mathcal{C}_S = \{x\in H^{1,1}(S,\mathbb{R}) | x\cdot x>0\:\:\text{and}\:\: x\cdot d >0
\:\:\text{for all nodal classes $d$}\}
\end{equation*}
\end{Th} 

We introduce further terms that allow us to describe the K\"ahler cone more
algebraically.

\begin{Def}
\begin{enumerate}
    \item Let $\ell_x\in \mathbb{P}(L_{\mathbb{C}})$. We define the 
    \emph{root system of} $\ell_x$ as

    \begin{equation*}
    \triangle_x := \{ d\in L| d\cdot d=-2,\: x\cdot d=0 \}\:.
    \end{equation*}
    
    \item We define the \emph{K\"ahler chambers} of $\ell_x$ as the
    connected components of 
    
    \begin{equation*}
    \{z\in L_{\mathbb{R}} | z\cdot z>0,\: z\cdot x=0,\: z\cdot d\neq 0
    \:\:\forall d\in\triangle_x\}\:. 
    \end{equation*}
\end{enumerate}    
\end{Def}

\begin{Th}
The subgroup of $\text{Aut}(L)$ that preserves $\ell_x$ 
acts transitively on the set of all K\"ahler chambers of $\ell_x$. The image 
$\phi_{\mathbb{R}}(\mathcal{C}_S)$ of the K\"ahler cone of a marked K3 
surface $(S,\phi)$ with period point $x$ is one of the K\"ahler chambers 
of $\ell_x$.    
\end{Th}

Finally, we arrive at the description of $\Kthree^m$. 

\begin{Def}
We define the \emph{augmented period domain} as

\begin{equation*}
\widetilde{\Omega} = \{(\ell_x,C) |\ell_x\in\Omega,\:\text{$C\subseteq L_{\mathbb{R}}$ 
is a K\"ahler chamber of $\ell_x$}\} 
\end{equation*}

and the \emph{augmented period map}\: $\widetilde{p}:\Kthree^m \rightarrow 
\widetilde{\Omega}$ by

\begin{equation*}
\widetilde{p}(S,\phi) := (p(S,\phi),\phi_{\mathbb{R}}(\mathcal{C}_S))\:.
\end{equation*}  
\end{Def}

\begin{Th}
The augmented period map $\widetilde{p}:\Kthree^m \rightarrow 
\widetilde{\Omega}$ is bijective. 
\end{Th}

In order to determine a Ricci-flat K\"ahler metric, we need to specify a 
K\"ahler class and not only the complex structure. Therefore, we need
another moduli space that takes this additional information into account. 

\begin{Def}
\begin{enumerate}
    \item A \emph{marked pair} is a pair of a marked K3 surface $(S,\phi)$
    and a K\"ahler class $y\in H^{1,1}(S,\mathbb{R})$. We usually write a marked pair as 
    $(S,\phi,y)$.
    
    \item Two marked pairs $(S,\phi,y)$ and $(S',\phi',y')$
    are called \emph{isomorphic} if there exists a biholomorphic map 
    $f:S\rightarrow S'$ that satisfies $\phi\circ f^{\ast} = \phi'$ and 
    $f^{\ast}y' = y$. 
    
    \item The \emph{moduli space of marked pairs $\Kthree^{mp}$}
    is the set of all marked pairs modulo isomorphisms.
\end{enumerate}
\end{Def}

Moreover, we define the following two sets:  

\begin{equation*}
\begin{array}{rcl}
K\Omega & := & \{ (\ell_x,y)\in \Omega \times L_{\mathbb{R}} | 
x\cdot y= 0, y\cdot y > 0 \} \\
K\Omega^0 & := & \{ (\ell_x,y)\in K\Omega | y\cdot d\neq 0 \:\: \forall d\in L
\:\:\text{with}\:\: d^2 = -2, x\cdot d=0 \} 
\end{array}
\end{equation*}

and the \emph{refined period map} 

\begin{eqnarray*}
p':\Kthree^{mp} & \rightarrow & \Omega \times L_{\mathbb{R}} \\
p'(S,\phi,y) & := & (p(S,\phi),\phi_{\mathbb{R}}(y))
\end{eqnarray*}

\begin{Th}
\label{TorelliRefined}
$p'$ takes its values in $K\Omega^0$. Moreover, it is a bijection 
between $\Kthree^{mp}$ and $K\Omega^0$. As a consequence, 
$\Kthree^{mp}$ is a real analytic Hausdorff ma\-nifold of
dimension $60$. 
\end{Th}

Finally, we describe the moduli space of all hyper-K\"ahler structures on K3 
surfaces. A hyper-K\"ahler structure on a marked K3 surface is a tuple
$(S,\phi,g,\omega_I,\omega_J,\omega_K)$, where $g$ is the hyper-K\"ahler 
metric and $\omega_I$, $\omega_J$ and $\omega_K$ are the K\"ahler
forms with respect to the complex structures $I$, $J$ and $K$ that satisfy
$IJK=-1$. The forms $\omega_I$, $\omega_J$ and $\omega_K$ determine
the metric and an orientation that makes $\omega_I\wedge\omega_J$ 
positive. Moreover, the cohomology classes of  $\omega_I$, $\omega_J$
and $\omega_K$ already determine $g$. $\phi_{\mathbb{C}}([\omega_J] + i[\omega_K])$
yields a period point of a K3 surface and $\phi_{\mathbb{R}}([\omega_I])$ determines a
K\"ahler chamber. This information determines the complex structure $I$
and then the K\"ahler class $[\omega_I]$ determines the hyper-K\"ahler
metric. Conversely, $g$ alone does only determine the span of $\omega_I$,
$\omega_J$ and $\omega_K$, but yields no basis of that space. The above 
observations motivate the following lemma on  isometries of K3 surfaces that 
will be useful later on. 

\begin{Le}
\label{IsomLem}
Let $S_j$ with $j\in\{1,2\}$ be K3 surfaces together with hyper-K\"ahler metrics $g_j$ 
and K\"ahler forms $\omega^{(j)}_I$, $\omega^{(j)}_J$ and $\omega^{(j)}_K$. Moreover,
let $V_j\subset H^2(S_j,\mathbb{R})$ be the subspace that is spanned by 
$[\omega^{(j)}_I]$, $[\omega^{(j)}_J]$ and $[\omega^{(j)}_K]$.  

\begin{enumerate}
    \item Let $f:S_1\rightarrow S_2$ be an isometry. The pull-back
    $f^{\ast}:H^2(S_2,\mathbb{Z}) \rightarrow H^2(S_1,\mathbb{Z})$ is
    a lattice isometry. Its $\mathbb{R}$-linear extension maps $V_2$
    to $V_1$.

    \item \label{Isom} Let $\psi:H^2(S_1,\mathbb{Z}) \rightarrow 
    H^2(S_2,\mathbb{Z})$ be a lattice isometry such that 
    $\psi_\mathbb{R}(V_1)=V_2$. Moreover, $\psi_\mathbb{R}$ shall
    map the positive cone of $S_1$ to the positive cone of $S_2$. Then 
    there exists an isometry $f:S_2\rightarrow S_1$ such that $f^{\ast} = \psi$. 

    \item Let $f:S\rightarrow S$ be an isometry that acts as the identity
    on $H^2(S,\mathbb{Z})$. Then, $f$ itself is the identity map. As a 
    consequence, the isometry from \ref{Isom}. is unique.  
\end{enumerate}
\end{Le}

\begin{proof}
The first claim is obvious and the third one follows from Proposition 11.3 in Chapter
VIII in \cite{BHPV}. The second claim is a consequence of the Torelli
theorem. More precisely, the fact that $\omega^{(j)}_J + i\omega^{(j)}_K$
is a $(2,0)$-form determines a splitting of $H^2(S_j,\mathbb{C})$ into
$H^{2,0}(S_j) \oplus H^{1,1}(S_j) \oplus H^{0,2}(S_j)$. Since $\psi_{\mathbb{R}}$
preserves the positive cone, it follows from our explicit description of
$\mathcal{C}_{S_j}$ that it preserves the K\"ahler cone, too. The Torelli theorem
thus yields a biholomorphic map $f:S_2\rightarrow S_1$ with $f^{\ast} = \psi$.
The triples $([\omega^{(2)}_I],[\omega^{(2)}_J],[\omega^{(2)}_K])$
and  $(\psi([\omega^{(1)}_I]),\psi([\omega^{(1)}_J]),\psi([\omega^{(1)}_K]))$
yield unique hyper-K\"ahler metrics on $S_2$. Since both triples span the same
subspace, these metrics are the same and we have $f^{\ast} g_1 = g_2$. 
\end{proof}

\begin{Rem}
If we had omitted the condition that $\psi_{\mathbb{R}}$ preserves the positive
cone, the second part of our lemma would have been slightly more complicated.
In that situation $\psi:=-\text{Id}_{H^2(S,\mathbb{Z})}$ would satisfy all conditions 
from the lemma. The corresponding isometry $f:S\rightarrow S$ would be
the identity map, but it would have to be interpreted as an anti-holomorphic map
between $(S,I)$ and $(S,-I)$. 
\end{Rem}

Finally, we describe the moduli space $\Kthree^{hk}$ of all marked 
hyper-K\"ahler structures $(S,\phi,g,\omega_I,$ $\omega_J,\omega_K)$. 
As a consequence of Theorem \ref{TorelliRefined} and Lemma \ref{IsomLem}
(see also  \cite[p. 161]{Joyce}), it follows that $\Kthree^{hk}$ is diffeomorphic
to the \emph{hyper-K\"ahler period domain} 

\begin{equation*}
\begin{aligned}
\Omega^{hk}:=\{ & (x,y,z)\in L_{\mathbb{R}}^3 | 
x^2 = y^2 = z^2>0,\: x\cdot y = x\cdot z = y\cdot z =0, \\
& \!\! \not\exists\: d\in L \:\:\text{with}\:\: d^2 = -2 \:\:\: \text{and} \:\:\:
x\cdot d = y\cdot d = z\cdot d=0    
\}\:. \\
\end{aligned}
\end{equation*}

\section{Singular K3 surfaces} In this section, we discuss singular K3 
surfaces and their relation to smooth ones. The results that
we present here were originally proven in \cite{Anderson1,Anderson2,Koba}.
A short overview can also be found in \cite[p.161 - 162]{Joyce}.

Let $S$ be a K3 surface and let $w\in H^2(S,\mathbb{Z})$ be a class 
with $w^2=-2$ that represents a submanifold $Z$ of $S$. We do not
assume that $w\in H^{1,1}(S)$ and thus $Z$ is not necessarily a
divisor. $S$ shall carry a hyper-K\"ahler structure $(g,\omega_I,\omega_J,
\omega_K)$. It can be shown that $Z$ can be chosen as a sphere that is
minimal with respect to $g$. Its area $A$ is given by 

\begin{equation*}
A^2 = \left([\omega_I]\cdot w\right)^2 + \left([\omega_J]\cdot w\right)^2
+ \left([\omega_K]\cdot w\right)^2
\end{equation*}

We choose a marking $\phi$ of $S$. If we move within the hyper-K\"ahler
period domain towards a triple $(x,y,z)\in L_{\mathbb{R}}^3$ with

\begin{equation*}
x\cdot \phi(w) = y\cdot \phi(w) = z\cdot \phi(w) = 0\:, 
\end{equation*}

the volume of the sphere shrinks to zero. In other words, we obtain a singularity. This
is in fact the geometric meaning of the condition in the definition of $\Omega^{hk}$ 
that there shall be no $d\in L$ with $d^2=-2$ and $x\cdot d = y\cdot d = z\cdot d =
0$. We assume that there is exactly one $d\in L$ with this property. In this situation,
we obtain the singularity by collapsing a single sphere with self-intersection $-2$
to a point. Since this is the reversal of blowing up an $A_1$-singularity, the K3 surface
has an $A_1$-singularity at a single point. Next, we assume that there exists an
arbitrary number of $d$s with $d^2=-2$ and $x\cdot d = y\cdot d = z\cdot d = 0$.
We define the set 

\begin{equation*}
\widetilde{\Omega}^{hk} := \{(x,y,z)\in L_{\mathbb{R}}^3 | 
x^2 = y^2 = z^2>0\:, x\cdot y = x\cdot z = y\cdot z=0 \}  
\end{equation*} 

and for any $\alpha=(x,y,z)\in \widetilde{\Omega}^{hk}$ we define

\begin{equation}
\label{SingSet}
\mathcal{D}_{\alpha}:= \{ d\in L | d^2=-2,\: x\cdot d = y\cdot d = 
z\cdot d = 0 \}\:.
\end{equation} 

By joining $d_1,d_2\in\mathcal{D}_{\alpha}$ with $d_1\neq d_2$ by
$d_1\cdot d_2$ edges, we obtain a graph $G$. This graph is the disjoint
union of simply laced Dynkin diagrams. As the hyper-K\"ahler structure
approaches $\alpha$, a set of 2-spheres whose intersection numbers are 
given by $d_i\cdot d_j$ collapses, which means that the Dynkin diagrams describe 
the type of the singularities. For example, if $G$ consists of one Dynkin diagram of 
type $E_8$ and $2$ isolated nodes, the singularities of the K3 surface are at $3$ 
different points. At one of them we have a singularity of type $E_8$ and at the other 
two ones we have $A_1$-singularities. Our considerations show that the singular
and the smooth marked K3 surfaces with a hyper-K\"ahler structure can be combined
into a larger moduli space that is diffeomorphic to $\widetilde{\Omega}^{hk}$.

\section{Non-symplectic involutions} 
\label{NonsympSection}
In this section, we introduce the most important results about non-symplectic
involutions. These results were proven by Nikulin \cite{Nikulin1,Nikulin2,Nikulin3} 
and are also summed up in \cite{ArSaTa,KoLe}. Moreover, we define a 
class of non-symplectic involutions that are well suited for explicit calculations
and we classify them.

\begin{Def}
Let $S$ be a K3 surface. A \emph{non-symplectic involution} is a biholomorphic 
map $\rho:S\rightarrow S$ such that

\begin{enumerate}
    \item $\rho^2 = \text{Id}$, but $\rho\neq \text{Id}$.
  
    \item The pull-back $\rho^{\ast}:H^{2,0}(S) \rightarrow H^{2,0}(S)$ is not 
    the identity map, or equivalently  $\rho^{\ast}(\omega_J +i\omega_K) = 
    - (\omega_J +i\omega_K)$.
\end{enumerate} 
\end{Def}

From now on, let $S$ be a K3 surface and $\rho:S\rightarrow S$ be a 
non-symplectic involution. We define the \emph{fixed lattice of $\rho$}
by 

\begin{equation*}
L^{\rho} := \{x\in H^2(S,\mathbb{Z}) | \rho^{\ast}x = x \}\:. 
\end{equation*}

$L^{\rho}$ is a primitive sublattice of $H^2(S,\mathbb{Z})$. Since
$\rho^{\ast}$ acts as $-1$ on $H^{2,0}(S)$ and $H^{0,2}(S)$,
$L^{\rho}$ is a sublattice of the Picard lattice.  A K3 surface with a 
non-symplectic involution admits an integral K\"ahler class $x$ and is thus 
algebraic by the Kodaira embedding theorem. Moreover, it admits an
integral $\rho$-invariant K\"ahler class since $x + \rho^{\ast} x$ is 
$\rho$-invariant. 

We choose a marking $\phi:H^2(S,\mathbb{Z}) \rightarrow L$
and abbreviate $\phi(L^{\rho})$ by $L^{\rho}$. It can be shown that
$L^{\rho}$ is a non-degenerate sublattice of $L$ with
signature $(1,t)$. A lattice with that kind of signature is called 
\emph{hyperbolic}. The rank $r=1+t$ is an invariant of
$L^{\rho}$. $L^{\rho}$ is \emph{2-elementary} which means that
$L^{\rho\ast}/L^{\rho}$ is isomorphic to a group of type 
$\mathbb{Z}_2^a$. The number $a\in\mathbb{N}_0$
is a second invariant of $L^{\rho}$. We define a third invariant
$\delta$ by

\begin{equation*}
\delta :=
\begin{cases}
0 & \text{if $x^2\in\mathbb{Z}$ for all $x\in L^{\rho\ast}$} \\
1 & otherwise \\
\end{cases}
\end{equation*} 

\begin{Th} 
\label{radelta-Theorem}
(Theorem 4.3.2 in \cite{Nikulin3})
For each triple $(r,a,\delta)\in \mathbb{N}_0\times\mathbb{N}_0
\times \{0,1\} $ there is up to isometries at most one even, 
hyperbolic, 2-elementary lattice with invariants $(r,a,\delta)$. 
\end{Th}

Let $N$ be a hyperbolic lattice such that there exists a primitive 
embedding of $N$ into $L$. We assume that $N^{\ast}/N$ is 
2-elementary and that $N\subset L$ contains a K\"ahler
class. Then there exists an involution $\rho_N$ of $L$ with fixed lattice 
$N$. $\rho_N$ acts as $-1$ on $N^{\perp}_\mathbb{R}
\subseteq L_{\mathbb{R}}$ and $N^{\perp}_{\mathbb{R}}$
contains a positive plane $P$ with an orthonormal basis $(x,y)$. 
The surjectivity of the period map and the Torelli
theorem guarantee that there exists a K3 surface $S$ together
with a non-symplectic involution $\rho$ such that $\rho^{\ast} =
\rho_N$ and $H^2(S,\mathbb{R}) \cap (H^{2,0}(S) \oplus
H^{0,2}(S)) = P$. The period point of that K3 surface is the
complex line that is spanned by $x+iy$.  

There is up to isometries of $L$ at most one primitive embedding 
of a lattice with invariants $(r,a,\delta)$ into $L$ and it follows that the 
deformation classes of K3 surfaces with a non-symplectic involution
can be classified in terms of triples $(r,a,\delta)$. Nikulin 
\cite{Nikulin3} has shown that there exist $75$ possible triples
that satisfy

\begin{equation*}
1\leq r\leq 20\:,\quad 0\leq a\leq 11\quad\text{and}\quad
r-a\geq 0\:.
\end{equation*}

A figure with a graphical representation of all possible values
of $(r,a,\delta)$ can be found in \cite{KoLe, Nikulin3}. 
Next, we describe the moduli space of all K3 surfaces with a 
non-symplectic involution whose fixed lattice is of a given isomorphism 
type. In order to do this, we need the following concept.

\begin{Def} (cf. Dolgachev \cite{Do2}) 
Let $N$ be a hyperbolic lattice that is primitively embedded
into $L$. 

\begin{enumerate}
    \item A \emph{marked ample $N$-polarised K3 surface} is a K3 
    surface $S$ together with a marking $\phi:H^2(S,\mathbb{Z}) 
    \rightarrow L$ such that $\phi^{-1}(N)$ is a sublattice of the Picard 
    lattice. Moreover, $\phi^{-1}(N)$ shall contain an integral ample 
    class, which is since $S$ is a compact K\"ahler manifold, the
    same as an integral K\"ahler class. 
    
    \item Two marked ample $N$-polarised K3 surfaces
    $(S,\phi)$ and $(S',\phi')$ are called isomorphic if there exists
    a biholomorphic map $f:S\rightarrow S'$ such that $\phi' =
    \phi \circ f^{\ast}$. 
    
    \item We denote the moduli space that consists of all marked
    ample $N$-polarised K3 surfaces modulo isomorphisms by 
    $\Kthree^m(N)$. 
\end{enumerate}     
\end{Def}

We denote the lattice with invariants $(r,a,\delta)$ by $L(r,a,\delta)$.
The moduli space of all marked K3 surfaces with a non-symplectic
involution whose fixed lattice is isomorphic to $L(r,a,\delta)$ 
is the same as $\Kthree^m(L(r,a,\delta))$, which we abbreviate by 
$\Kthree^m(r,a,\delta)$. We remark that this moduli space is the
same as the moduli space $\Kthree'(r,a,\delta)$ in \cite{KoLe}.
There is a nice explicit description of $\Kthree^m(r,a,\delta)$.

\begin{Th} (Corollary 3.2 in \cite{Do2})
Let $N$ be a hyperbolic lattice that can be primitively embedded 
into $L$. We denote the orthogonal complement of
$N$ in $L$ by $M$ and define the following sets:

\begin{equation*}
\begin{array}{rcl}
\Omega_N & := & \{\ell_x \in \mathbb{P}(M_{\mathbb{C}}) | x\cdot x = 0,
x\cdot\overline{x} > 0 \} \\
\triangle(M) & := & \{d\in M | d^2 = -2 \} \\
H_d & := & \{\ell_z\in \mathbb{P}(M_{\mathbb{C}}) | z\cdot d = 0\} \\
\Omega'_N & := & \Omega_N \setminus \bigcup_{d\in\triangle(M)}
(H_d \cap \Omega_N) \\  
\end{array}
\end{equation*} 

$\Kthree^m(N)$ is isomorphic to $\Omega'_N$ and the isomorphism
is given by the period map.
\end{Th}

\begin{Rem}
Since $N$ contains a K\"ahler class, any element of $H_d\cap \Omega_N$ would
correspond to a K3 surface $S$ with the property that $d$ is orthogonal
to the real and imaginary part of the $(2,0)$-form and to a K\"ahler
class. In other words, $S$ would carry a singular hyper-K\"ahler metric. 
This is the reason why we have to remove the set $\bigcup_{d\in\triangle(M)}
(H_d\cap\Omega_N)$ from $\Omega_N$.
\end{Rem}

The topology of the fixed locus $S^{\rho}:=\{x\in S | \rho(x)=x \}$ of a
non-symplectic involution $\rho$ can be described in terms of the invariants $r$
and $a$. 

\begin{Th} \label{FixedLocusTheorem}
(cf. \cite{KoLe,Nikulin3})
Let $\rho:S\rightarrow S$ be a non-symplectic involution of a K3 surface
and let $(r,a,\delta)$ be the invariants of the fixed lattice. The fixed locus 
$S^{\rho}$ of $\rho$ is a disjoint union of complex curves.

\begin{enumerate}
    \item If $(r,a,\delta)=(10,10,0)$, $S^{\rho}$ is empty.
    
    \item If $(r,a,\delta)=(10,8,0)$, $S^{\rho}$ is the disjoint union
    of two elliptic curves.
    
    \item In the remaining cases, we have
    
    \begin{equation*}
    S^{\rho} = C_g \cup E_1 \cup \ldots \cup E_k\:,
    \end{equation*}  
    
    where $C_g$ is a curve of genus $g=\tfrac{22 - r - a}{2}$ and
    the $E_i$ are $k=\tfrac{r-a}{2}$ curves that are biholomorphic  
    to $\mathbb{CP}^1$, i.e. they are rational curves.   
\end{enumerate}
\end{Th} 

We define a class of non-symplectic involutions whose action on 
$L$ has a very simple matrix representation. In order to do this, 
we have to fix a basis of $L$. We write

\begin{equation*}
L = H_1 \oplus H_2 \oplus H_3 \oplus (-E_8)_1 \oplus (-E_8)_2
\end{equation*}

in order to distinguish between the different summands. We 
choose a basis $(u_1^i,u_2^i)$ of each $H_i$ such that 

\begin{equation*}
u_1^i\cdot u_1^i = u_2^i\cdot u_2^i = 0\:,
\quad
u_1^i\cdot u_2^i = 1\:. 
\end{equation*}

Moreover, $(v_1^i,\ldots,v_8^i)$ shall be a basis of 
$(-E_8)_i$ such that the bilinear form has the matrix 
representation 

\begin{equation*}
\left(\,\begin{array}{cccccccc}
-2 & 0 & 1 & 0 & 0 & 0 & 0 & 0\\
0 & -2 & 0 & 1 & 0 & 0 & 0 & 0\\
1 & 0 & -2 & 1 & 0 & 0 & 0 & 0\\ 
0 & 1 & 1 & -2 & 1 & 0 & 0 & 0\\
0 & 0 & 0 & 1 & -2 & 1 & 0 & 0\\
0 & 0 & 0 & 0 & 1 & -2 & 1 & 0\\
0 & 0 & 0 & 0 & 0 & 1 & -2 & 1\\
0 & 0 & 0 & 0 & 0 & 0 & 1 & -2\\
\end{array}\,\right)
\end{equation*}

We call 

\[
(w_1,\ldots,w_{22}) = (u_1^1,u_2^1,u_1^2,u_2^2,u_1^3,u_2^3,
v_1^1,\ldots,v_8^1,v_1^2,\ldots,v_8^2)
\]

the \emph{standard basis of $L$}. With help of this basis, 
we are able to define our class of non-symplectic involutions.

\begin{Def}
Let $S$ be a K3 surface and let $\rho:S\rightarrow S$  be
a non-symplectic involution. We call $\rho$ a \emph{simple
non-symplectic involution} if there exists a marking 
$\phi:H^2(S,\mathbb{Z})\rightarrow L$ such that for all
$i\in\{1,\ldots,22\}$ there exists a $j\in\{1,\ldots,22\}$
with $\rho(w_i)= \pm w_j$, where 
$\phi\circ\rho^{\ast}\circ\phi^{-1}$ is abbreviated
by $\rho$.  
\end{Def} 

Let $\rho$ be a simple non-symplectic involution. Since 
$\rho:L\rightarrow L$ is a lattice isometry and we have
$\rho(w_i)=\pm w_j$, $\rho$ maps any of the sublattices 
$H_k\subseteq L$ to an $H_l$.  There are four possibilities
for the value of $\rho(u_1^k)$ and of $\rho(u_2^k)$. We check 
for each combination if $\rho|_{H_k}: H_k\rightarrow H_l$ is a 
lattice isometry and see that $\rho|_{H_k}$ is given by one of the 
following maps:

\begin{enumerate}
    \item $\rho|_{H_k}(u_1^k) = u_1^l$, $\rho|_{H_k}(u_2^k) = u_2^l$
    
    \item $\rho|_{H_k}(u_1^k) = -u_1^l$, $\rho|_{H_k}(u_2^k) = -u_2^l$
    
    \item $\rho|_{H_k}(u_1^k) = u_2^l$, $\rho|_{H_k}(u_2^k) = u_1^l$
    
    \item $\rho|_{H_k}(u_1^k) = -u_2^l$, $\rho|_{H_k}(u_2^k) = -u_1^l$
\end{enumerate}

Since $\rho$ is non-symplectic, its fixed lattice is hyperbolic. 
Therefore, $\rho|_{3H}:3H\rightarrow 3H$ has to preserve exactly 
one positive vector. By enumerating all possibilities for $\rho|_{3H}$ 
with this property and comparing the invariants of the fixed lattices, 
we can conclude that $\rho|_{3H}$ is up to conjugation one of the 
following maps $\rho_1^i: 3H\rightarrow 3H$ with $i=1,\ldots,7$:  

\bigskip

\begin{center}
\label{TableRho1}
\begin{longtable}{l|l|l|l|l|l}

$i$ & Matrix representation & Basis of the fixed lattice & $r$ & $a$ & $\delta$ \\

\hline \hline 

&&&& \\

$1$ & 
$\left(\,
\begin{array}{cccccc}
\hhline{--~~~~}
\multicolumn{1}{|c}{1} & \multicolumn{1}{c|}{0} & & & & \\
\multicolumn{1}{|c}{0} & \multicolumn{1}{c|}{1} & & & & \\
\hhline{----~~}
& & \multicolumn{1}{|c}{-1} & \multicolumn{1}{c|}{0} & & \\
& & \multicolumn{1}{|c}{0} & \multicolumn{1}{c|}{-1} & & \\
\hhline{~~----}
& & & & \multicolumn{1}{|c}{-1} & \multicolumn{1}{c|}{0} \\
& & & & \multicolumn{1}{|c}{0} & \multicolumn{1}{c|}{-1} \\
\hhline{~~~~--}
\end{array}
\,\right)$ & $(u_1^1,u_2^1)$ & $2$ & $0$ & $0$ \\

&&&& \\
\hline
&&&& \\\nopagebreak[4]

$2$ & 
$\left(\,
\begin{array}{cccccc}
\hhline{--~~~~}
\multicolumn{1}{|c}{1} & \multicolumn{1}{c|}{0} & & & & \\
\multicolumn{1}{|c}{0} & \multicolumn{1}{c|}{1} & & & & \\
\hhline{----~~}
& & \multicolumn{1}{|c}{0} & \multicolumn{1}{c|}{-1} & & \\
& & \multicolumn{1}{|c}{-1} & \multicolumn{1}{c|}{0} & & \\
\hhline{~~----}
& & & & \multicolumn{1}{|c}{-1} & \multicolumn{1}{c|}{0} \\
& & & & \multicolumn{1}{|c}{0} & \multicolumn{1}{c|}{-1} \\
\hhline{~~~~--}
\end{array}
\,\right)$ & $(u_1^1,u_2^1,u_1^2-u_2^2)$ & $3$ & $1$ & $1$ \\

&&&& \\
\hline
&&&& \\\nopagebreak[4]

$3$ & 
$\left(\,
\begin{array}{cccccc}
\hhline{--~~~~}
\multicolumn{1}{|c}{1} & \multicolumn{1}{c|}{0} & & & & \\
\multicolumn{1}{|c}{0} & \multicolumn{1}{c|}{1} & & & & \\
\hhline{----~~}
& & \multicolumn{1}{|c}{0} & \multicolumn{1}{c|}{-1} & & \\
& & \multicolumn{1}{|c}{-1} & \multicolumn{1}{c|}{0} & & \\
\hhline{~~----}
& & & & \multicolumn{1}{|c}{0} & \multicolumn{1}{c|}{-1} \\
& & & & \multicolumn{1}{|c}{-1} & \multicolumn{1}{c|}{0} \\
\hhline{~~~~--}
\end{array}
\,\right)$ & $(u_1^1,u_2^1,u_1^2-u_2^2,u_1^3-u_2^3)$ & $4$ & $2$ & $1$ \\

&&&& \\
\hline
&&&& \\\nopagebreak[4]

$4$ & 
$\left(\,
\begin{array}{cccccc}
\hhline{--~~~~}
\multicolumn{1}{|c}{0} & \multicolumn{1}{c|}{1} & & & & \\
\multicolumn{1}{|c}{1} & \multicolumn{1}{c|}{0} & & & & \\
\hhline{----~~}
& & \multicolumn{1}{|c}{-1} & \multicolumn{1}{c|}{0} & & \\
& & \multicolumn{1}{|c}{0} & \multicolumn{1}{c|}{-1} & & \\
\hhline{~~----}
& & & & \multicolumn{1}{|c}{-1} & \multicolumn{1}{c|}{0} \\
& & & & \multicolumn{1}{|c}{0} & \multicolumn{1}{c|}{-1} \\
\hhline{~~~~--}
\end{array}
\,\right)$ & $(u_1^1+u_2^1)$ & $1$ & $1$ & $1$ \\

&&&& \\
\hline
&&&& \\\nopagebreak[4]

$5$ & 
$\left(\,
\begin{array}{cccccc}
\hhline{--~~~~}
\multicolumn{1}{|c}{0} & \multicolumn{1}{c|}{1} & & & & \\
\multicolumn{1}{|c}{1} & \multicolumn{1}{c|}{0} & & & & \\
\hhline{----~~}
& & \multicolumn{1}{|c}{0} & \multicolumn{1}{c|}{-1} & & \\
& & \multicolumn{1}{|c}{-1} & \multicolumn{1}{c|}{0} & & \\
\hhline{~~----}
& & & & \multicolumn{1}{|c}{-1} & \multicolumn{1}{c|}{0} \\
& & & & \multicolumn{1}{|c}{0} & \multicolumn{1}{c|}{-1} \\
\hhline{~~~~--}
\end{array}
\,\right)$ & $(u_1^1+u_2^1,u_1^2-u_2^2)$ & $2$ & $2$ & $1$ \\ 

&&&& \\
\hline 
&&&& \\\nopagebreak[4]

$6$ & 
$\left(\,
\begin{array}{cccccc}
\hhline{--~~~~}
\multicolumn{1}{|c}{0} & \multicolumn{1}{c|}{1} & & & & \\
\multicolumn{1}{|c}{1} & \multicolumn{1}{c|}{0} & & & & \\
\hhline{----~~}
& & \multicolumn{1}{|c}{0} & \multicolumn{1}{c|}{-1} & & \\
& & \multicolumn{1}{|c}{-1} & \multicolumn{1}{c|}{0} & & \\
\hhline{~~----}
& & & & \multicolumn{1}{|c}{0} & \multicolumn{1}{c|}{-1} \\
& & & & \multicolumn{1}{|c}{-1} & \multicolumn{1}{c|}{0} \\
\hhline{~~~~--}
\end{array}
\,\right)$ & $(u_1^1+u_2^1,u_1^2-u_2^2,u_1^3-u_2^3)$ & $3$ & $3$ & $1$ \\

&&&& \\
\hline
&&&& \\\nopagebreak[4]

$7$ & 
$\left(\,
\begin{array}{cccccc}
\hhline{----~~}
\multicolumn{1}{|c}{0} & 0 & 1 & \multicolumn{1}{c|}{0} & & \\
\multicolumn{1}{|c}{0} & 0 & 0 & \multicolumn{1}{c|}{1} & & \\
\multicolumn{1}{|c}{1} & 0 & 0 & \multicolumn{1}{c|}{0} & & \\
\multicolumn{1}{|c}{0} & 1 & 0 & \multicolumn{1}{c|}{0} & & \\
\hhline{------}
& & & & \multicolumn{1}{|c}{-1} & \multicolumn{1}{c|}{0} \\
& & & & \multicolumn{1}{|c}{0} & \multicolumn{1}{c|}{-1} \\
\hhline{~~~~--}
\end{array}
\,\right)$ & $(u_1^1+u_1^2,u_2^1 + u_2^2)$ & $2$ & $2$ & $0$ \\
\end{longtable} 
\end{center}

We study the restriction of $\rho$ to $2(-E_8)$. Let $i\in\{7,\ldots,22\}$. 
If $\rho(w_i)=w_j$ with $i\neq j$, we have $\rho(w_j)=w_i$
since $\rho$ is an involution. If $\rho(w_i)=-w_j$ with $i\neq j$, 
we have $\rho(w_j)=-w_i$ for the same reason. Therefore, there
exists a permutation $\sigma$ of $\{7,\ldots,22\}$ such that 
the basis $(w'_1,\ldots,w'_{16}) := (w_{\sigma(7)},\ldots,w_{\sigma(22)})$ 
satisfies:

\begin{enumerate}
    \item $\rho(w'_i)=w'_i$ for $i\in\{1,\ldots,k_1\}$,
    
    \item $\rho(w'_i)=-w'_i$ for $i\in\{k_1+1,\ldots,k_2\}$,
    
    \item $\rho(w'_{2i-1})= w'_{2i}$ and $\rho(w'_{2i})= w'_{2i-1}$
    for $i\in\{\frac{k_2}{2}+1,\ldots,k_3\}$ and
    
    \item $\rho(w'_{2i-1})= -w'_{2i}$ and $\rho(w'_{2i})= -w'_{2i-1}$
    for $i\in\{k_3+1,\ldots,8\}$.
\end{enumerate}

for suitable $k_1,k_2,k_3\in\mathbb{N}_0$. Let $i\in\{k_1+1,\ldots,k_2\}$, which
means that $\rho(w'_i)=-w'_i$. The number $i$ corresponds to a node of one of the 
two Dynkin diagrams of type $E_8$. Let $j$ be a node that is connected to $i$ by 
an edge. The restriction of the bilinear form to $span(w'_i,w'_j)$ is given 
by

\[
\left(\,
\begin{array}{cc}
-2 & 1 \\
1 & -2 \\
\end{array}
\,\right)
\]

If $\rho(w'_j)= w'_j$, $\rho$ does not preserve the bilinear form. Therefore, we have
$\rho(w'_j)=\pm w'_k$ with $k\neq i,j$ or $\rho(w'_j)=- w'_j$. We assume that
$\rho(w'_j)=\pm w'_k$. Since $-w'_i\cdot \pm w'_k = \rho(w'_i)\cdot
\rho(w'_j) = w'_i\cdot w'_j = 1$ and all off-diagonal coefficients of the
Cartan matrix are positive, we have $\rho(w'_j)=- w'_k$ and $i$ and $k$ 
have to be connected by an edge. Analogously, we can conclude that any 
node that is connected to $j$ is mapped to a node that is connected to $k$. 
By repeating this argument, it follows that $\rho$ acts as a non-trivial graph 
automorphism on the diagram $E_8$ to which $i$ belongs. Since $E_8$ 
has no symmetries, this is impossible and we have $\rho(w'_j)=-w'_j$.
Again, we can repeat this argument and conclude that $\{k_1+1,\ldots,k_2\}$
consists of zero, one or both connected components of $2E_8$. 

Next, let $i\in\{\frac{k_2}{2}+1,\ldots,k_3\}$, which means that $w'_{2i-1}$ is 
mapped to another basis element $w'_{2i}$. By the same argument as above, 
we see that all nodes that are connected to $2i-1$ are mapped to nodes that 
are connected to $2i$. The restriction of $\rho$ to $span(w'_{k_2+1},
\ldots,w'_{2k_3})$ thus maps connected components of $2E_8$ to other 
connected components. It follows that either $\{\tfrac{k_2}{2}+1,\ldots,k_3\}$ 
is empty or $\rho$ interchanges both copies of $E_8$. Finally, let 
$i\in\{k_3+1,\ldots,8\}$. In this case, we have $\rho(w'_{2i-1})= -w'_{2i}$ 
and it follows that if $\{k_3+1,\ldots,8\}$ is not empty, the first $E_8$ is 
mapped to the second $E_8$ such that $v_k^1$ is mapped to $-v_k^2$. 
All in all, the restricted map $\rho|_{2(-E_8)}: 2(-E_8) \rightarrow 2(-E_8)$ 
is up to conjugation one of the involutions $\rho_2^j$ 
below. Let $x_1\in (-E_8)_1$ and $x_2\in (-E_8)_2$. We define 

\begin{equation}
\begin{array}{ll}
\rho_2^1(x_1,x_2):=(x_1,x_2)\:,\quad & 
\rho_2^2(x_1,x_2):=(-x_1,x_2)\:,\\
\rho_2^3(x_1,x_2):=(-x_1,-x_2)\:,\quad & 
\rho_2^4(x_1,x_2):=(x_2,x_1)\:. \\
\end{array} 
\end{equation}    

Moreover, any conjugate $\psi: 2(-E_8) \rightarrow 2(-E_8)$ of the 
$\rho_2^j$ that is still simple is given either by  

\begin{equation}
\label{rho2alt}
\psi(x_1,x_2)=(x_1,-x_2) 
\quad\text{or}\quad
\psi(x_1,x_2)=(-x_2,-x_1) 
\end{equation}

The fixed lattices and invariants of the $\rho_2^j$ can be found 
in the following table:

\begin{center}
\begin{tabular}{l|c|l|l|l}
$j$ & Fixed lattice & $r$ & $a$ & $\delta$ \\

\hline

$1$ & $2(-E_8)$ & $16$ & $0$ & $0$ \\

$2$ & $-E_8$ & $8$ & $0$ & $0$ \\

$3$ & $\{0\}$ & $0$ & $0$ & $0$ \\

$4$ & $-E_8(2)$ & $8$ & $8$ & $0$ \\
\end{tabular}
\end{center}

Any $\rho_1^i \oplus \rho_2^j$ with $1\leq i\leq 7$ and
$1\leq j\leq 4$ is an involution of $L$. Since the complement 
of the fixed lattice contains a positive plane, we can conclude with
help of the Torelli theorem or with Lemma \ref{IsomLem}
that these lattice involutions are pull-backs of non-symplectic
involutions. Finally, we compute the invariants of the
involutions that we have found. Let $K=K_1\oplus K_2$
be a direct sum of even, hyperbolic, 2-elementary lattices.
We denote the invariants of $K$ by $(r,a,\delta)$
and those of the $K_i$  by $(r_i,a_i,\delta_i)$. It is easy
to see that $r=r_1+r_2$, $a=a_1+a_2$ and
$\delta= \max\{\delta_1,\delta_2\}$.  Therefore, we
have proven the following theorem:

\begin{Th} \label{SimpleNonSymplThm}
Let $S$ be a K3 surface and let $\rho:S\rightarrow S$
be a non-symplectic involution. $\rho$ is simple if 
and only if its invariants $(r,a,\delta)$ can be found
in the table below. Moreover, the action of $\rho$
on the K3 lattice $L$ is conjugate to an involution
$\rho_1^i \oplus \rho_2^j$ that we have defined above.
The values of $i$ and $j$ that correspond to an
involution with invariants $(r,a,\delta)$ are also
included in the following table. 

\begin{center}
\begin{tabular}{llll}
\begin{tabular}{l|l}
$(i,j)$ & $(r,a,\delta)$ \\
\hline
$(1,1)$ & $(18,0,0)$ \\
$(1,2)$ & $(10,0,0)$ \\
$(1,3)$ & $(2,0,0)$ \\
$(1,4)$ & $(10,8,0)$ \\
$(2,1)$ & $(19,1,1)$ \\
$(2,2)$ & $(11,1,1)$ \\
$(2,3)$ & $(3,1,1)$ \\
$(2,4)$ & $(11,9,1)$ \\
$(3,1)$ & $(20,2,1)$ \\
$(3,2)$ &  $(12,2,1)$ \\
$(3,3)$ & $(4,2,1)$ \\
$(3,4)$ & $(12,10,1)$ \\
$(4,1)$ & $(17,1,1)$ \\
$(4,2)$ & $(9,1,1)$ \\
\end{tabular} & 
\begin{tabular}{l|l}
$(i,j)$ & $(r,a,\delta)$ \\
\hline
$(4,3)$ & $(1,1,1)$ \\
$(4,4)$ & $(9,9,1)$ \\
$(5,1)$ & $(18,2,1)$ \\
$(5,2)$ & $(10,2,1)$ \\
$(5,3)$ & $(2,2,1)$ \\
$(5,4)$& $(10,10,1)$ \\
$(6,1)$ & $(19,3,1)$ \\
$(6,2)$ & $(11,3,1)$ \\
$(6,3)$ & $(3,3,1)$ \\
$(6,4)$ & $(11,11,1)$ \\
$(7,1)$ & $(18,2,0)$ \\
$(7,2)$ & $(10,2,0)$ \\
$(7,3)$ & $(2,2,0)$ \\
$(7,4)$ & $(10,10,0)$ \\
\end{tabular} 
\end{tabular}
\end{center}
\end{Th}

\section{K3 surfaces with singularities and a non-symplectic involution}
\label{OneInvolution}
In this section, we study which kinds of ADE-singularities a K3 surface
with a non-symplectic involution may have. We focus
on the case where the non-symplectic involution is simple. 
Let $(S,\phi)$ be a marked K3 surface with a hyper-K\"ahler 
structure and a distinguished complex structure. Moreover, let $\ell_{x+iy}$ 
be its period point and let $z\in L$ be the image of the K\"ahler class with respect 
to $\phi$. We assume that $S$ admits a non-symplectic involution $\rho$ 
with invariants $(r,a,\delta)$ that leaves the metric invariant. This implies that 

\begin{equation}
\label{KaehlerRelations1}
\rho(x)=-x\:,\quad  \rho(y)=-y\:,\quad \rho(z)=z\:.
\end{equation}

We recall that the set 

\begin{equation*}
D:= \{d\in L | d^2=-2\:, x\cdot d=y\cdot d=z\cdot d=0 \}
\end{equation*}

is a root system that determines the number and type of the singular 
points. In order to study the possible singularities of $S$, we choose 
$x$, $y$ and $z$ in such a way that $D$ is large. $x$ and $y$ have to be 
positive elements in the orthogonal complement of the fixed lattice $L^{\rho}$.
Our description of the moduli space $\Kthree^m(r,a,\delta)$ guarantees
that any choice of $x,y \in L^{\rho\perp}$ with $x^2=y^2>0$ and
$x\cdot y=0$ yields a period point of a (possibly singular) K3 surface with 
a non-symplectic involution with invariants $(r,a,\delta)$. Moreover, we
can choose $z$ as an arbitrary element of $L^{\rho}$ with $z^2=x^2$
and $z\cdot x = z\cdot y = 0$. 

We assume that $\rho$ is a simple non-symplectic involution 
and that we have chosen the marking such that $\rho$ acts as
$\rho^i_1 \oplus \rho^j_2$ on $L$. Depending on $i$ we
choose $z\in L_{\mathbb{R}}$ as follows:

\[
z:=
\begin{cases}
u_1^1 + u_2^1 & \text{if}\:\:\: 1\leq i\leq 6, \\
u_1^1 + u_2^1 + u_1^2 + u_2^2 & \text{if}\:\:\: i=7. \\
\end{cases}
\]

If $i=7$, we have $z^2=4$ and we have $z^2=2$ otherwise. We choose 
$x$ and $y$ as:

\[
\begin{array}{l}
x:=
\begin{cases}
u_1^2 + u_2^2 & \text{if}\:\:\: 1\leq i\leq 6, \\
u_1^1 + u_2^1 - u_1^2 - u_2^2 & \text{if}\:\:\: i=7. \\
\end{cases} \\ \\
y:=
\begin{cases}
u_1^3 + u_2^3 & \text{if}\:\:\: 1\leq i\leq 6, \\
\sqrt{2}(u_1^3 + u_2^3) & \text{if}\:\:\: i=7. \\
\end{cases} \\
\end{array}
\]

$x$, $y$ and $z$ satisfy $x^2=y^2=z^2$ and the three vectors
are pairwise orthogonal. By a short calculation, we see that  
for any value of $i$ we have $z\in L^{\rho}$ and $x$ as well as 
$y$ is orthogonal to $L^{\rho}$. The orthogonal complement
of $span_{\mathbb{Z}}(x,y,z)$ is for all values of $i$ given by

\begin{equation}
\label{PicardMax}
span_{\mathbb{Z}}(u_1^1 - u_2^1,u_1^2 - u_2^2,
u_1^3 - u_2^3) \oplus (-E_8)_1 \oplus (-E_8)_2\:.
\end{equation}

A K3 surface $S$ with a hyper-K\"ahler structure that is determined by $x$, 
$y$ and $z$ thus has $3$ singular points of type $A_1$ and $2$ singular 
points of type $E_8$. Since $\rho$ acts on $x$, $y$ and $z$ as in equation 
(\ref{KaehlerRelations1}), there exists a non-symplectic involution of 
$S$ with fixed lattice $L^{\rho}$. All in all, we have proven the following 
theorem.  

\begin{Th}
\label{K3SingThm}
Let $(r,a,\delta)\in\mathbb{N}\times\mathbb{N}_0\times \{0,1\}$ be 
a triple such that there exists a K3 surface with a simple non-symplectic involution 
with invariants $(r,a,\delta)$. Then there exists a K3 surface which 
has $3$ singular points with $A_1$-singularities and $2$ singular points with 
$E_8$-singularities and carries a hyper-K\"ahler metric that is invariant with 
respect to a non-symplectic involution with the same values of $(r,a,\delta)$. 
\end{Th}

\begin{Rem}
The Picard lattice of the K3 surface from the above theorem is the direct sum
of the lattice (\ref{PicardMax}) and $span_{\mathbb{Z}}(z)$ and $S$ therefore 
has maximal Picard number. 
\end{Rem}

We search for K3 surfaces with a simple non-symplectic involution whose
singularities are of a different kind.  More precisely, let $G$ be a Dynkin diagram 
that can be obtained by deleting some nodes from the union of three Dynkin
diagrams of type $A_1$ and two of type $E_8$. We investigate if there exists 
a K3 surface with a simple non-symplectic involution whose singularities are
described by $G$. We denote the lattice (\ref{PicardMax}) by $N$ and fix a basis 

\begin{equation}
\label{PicardMaxBasis}
(\widetilde{w}_1,\ldots,\widetilde{w}_{19}):=(u_1^1-u_2^1,u_1^2-u_2^2,u_1^3-u_2^3,
v_1^1,\ldots,v_8^1,v_1^2,\ldots,v_8^2)
\end{equation}

of $N$. Let $S$ be a K3 surface with a simple non-symplectic involution $\rho$ whose
invariants are $(r,a,\delta)$. We choose a marking such that $\rho(w_i)=\pm w_j$. 
It is easy to see that $\rho(N)=N$ and that for any $i\in\{1,\ldots,19\}$ there
exists a $j$ such that $\rho(\widetilde{w}_i)=\widetilde{w}_j$. For the same reasons as 
in Section \ref{NonsympSection}, there exists a permutation $\sigma$  of 
$\{1,\ldots,19\}$ such that $(\widetilde{w}'_1,\ldots,\widetilde{w}'_{19}) := 
(\widetilde{w}_{\sigma(1)},\ldots,\widetilde{w}_{\sigma(19)})$ satisfies:

\begin{enumerate}
    \item $\rho(\widetilde{w}'_i)=\widetilde{w}'_i$ for $i\in\{1,\ldots,k_1\}$,
    
    \item $\rho(\widetilde{w}'_i)=-\widetilde{w}'_i$ for $i\in\{k_1+1,\ldots,k_2\}$,
    
    \item $\rho(\widetilde{w}'_{2i})= \widetilde{w}'_{2i+1}$ and 
    $\rho(\widetilde{w}'_{2i+1})= \widetilde{w}'_{2i}$ for $i\in\{\frac{k_2+1}{2},
    \ldots,k_3\}$ and
    
    \item $\rho(\widetilde{w}'_{2i})= -\widetilde{w}'_{2i+1}$ and 
    $\rho(\widetilde{w}'_{2i+1})= -\widetilde{w}'_{2i}$ for $i\in\{k_3+1,\ldots,9\}$.
\end{enumerate}

for suitable $k_1,k_2,k_3\in\mathbb{N}_0$. We choose four arbitrary subsets 
$M_1\subseteq \{1,\ldots,k_1\}$, $M_2\subseteq \{k_1+1,\ldots,k_2\}$, 
$M_3\subseteq \{\frac{k_2+1}{2},\ldots,k_3\}$ and $M_4\subseteq\{k_3+1,\ldots,9\}$.
Moreover, we choose for any $j\in M_i$ an $\alpha_{ij}\in\mathbb{R}$ such that the 
family

\[ 
(1,\alpha_{1\min{M_1}},\ldots,\alpha_{1\max{M_1}},\ldots,
\alpha_{4\min{M_4}},\ldots,\alpha_{4\max{M_4}})
\]

is $\mathbb{Q}$-linearly independent. We replace $x,y,z\in L_{\mathbb{R}}$
that we have defined in the proof of Theorem \ref{K3SingThm} by

\begin{equation}
\begin{array}{rcl}
x' & =& x + \sum_{j\in M_2} \alpha_{2j} \widetilde{w}'_j
+ \sum_{j\in M_3} \alpha_{3j} (\widetilde{w}'_{2j} - \widetilde{w}'_{2j+1})
+ \sum_{j\in M_4} \alpha_{4j} (\widetilde{w}'_{2j} + \widetilde{w}'_{2j+1}) \\
y' & = & \left(\frac{x'^2}{y^2}\right)^{\tfrac{1}{2}} y \\
z' & = & z +  \sum_{j\in M_1} \alpha_{1j} \widetilde{w}'_j 
\end{array}
\end{equation}

$x'$ and $y'$ are still in the $(-1)$-eigenspace of $\rho$ and 
$z'$ is still $\rho$-invariant. If the $\alpha_{ij}$ are sufficiently small,
$x'$, $y'$ and $z'$ are positive. We have

\begin{equation*}
x'^2 = x^2 - 2\sum_{j\in M_2} \alpha_{2j}^2  - 4\sum_{j\in M_3} \alpha_{3j}^2
- 4\sum_{j\in M_4} \alpha_{4j}^2 = y'^2\:, \quad
z'^2 = z^2 - 2\sum_{j\in M_1} \alpha_{ij}^2\:,
\end{equation*}

since $x$, $y$ and $z$ are orthogonal to $N$. It is possible to choose
the $\alpha_{ij}$ such that 

\begin{equation*}
2\sum_{j\in M_2} \alpha_{2j}^2 + 4\sum_{j\in M_3} \alpha_{3j}^2
+ 4\sum_{j\in M_4} \alpha_{4j}^2 = 2\sum_{j\in M_1} \alpha_{ij}^2
\end{equation*}

and thus we can assume that $x'^2=y'^2=z'^2>0$. Moreover, we have
$x'\cdot y'=x'\cdot z'=y'\cdot z'=0$. If $M_1=\emptyset$ or 
$M_2\cup M_3\cup M_4 = \emptyset$, we can define
$z'=\lambda z$ or $x'=\mu x$ and $y'= \mu y$ for appropriate 
$\lambda,\mu\in\mathbb{R}$ such that $x'^2=y'^2=z'^2$.
Therefore, we obtain a triple $(x',y',z')$ with the same properties as
above in that case.

All in all, the triple $(x',y',z')$ defines a new hyper-K\"ahler 
structure on $S$. Since $x'$ and $y'$ remain in the $(-1)$-eigenspace of 
$\rho$, $S$ admits a non-symplectic involution with the same fixed lattice
as before. Since $z'$ is $\rho$-invariant, $\rho$ is the pull-back of an isometry
with respect to the new hyper-K\"ahler metric. The set

\begin{equation*}
\begin{aligned}
D' & := \{ d\in L | d^2=-2, x'\cdot d= y'\cdot d = z'\cdot d = 0 \} \\
& = \left\{ \widetilde{w}'_i  \middle| i\notin M_1\cup M_2 \wedge \frac{i}{2}\notin M_3\cup M_4 \wedge \frac{i-1}{2}
\notin M_3\cup M_4 \right\}
\end{aligned}
\end{equation*}

is a root system that describes the number and type of the singular points
of the new hyper-K\"ahler metric. 

We interpret $D'$ geometrically. Any $\widetilde{w}'_i$ with $i\notin M_1\cup M_2$
corresponds to a sphere $S^2$ with vanishing area. The isometry $\rho:S \rightarrow
S$ maps such an $S^2$ to another $S^2$ with vanishing area. Since 
$\rho(\widetilde{w}'_i)=\pm \widetilde{w}'_i$, the sphere is mapped to itself and the 
sign determines if $\rho$ acts orientation-preserving on the sphere. Analogously, the 
$\widetilde{w}'_{2i}$ with  $i\notin M_3\cup M_4$ and the $\widetilde{w}'_{2i+1}$ with
$i\notin M_3\cup M_4$ correspond to sets of spheres with area $0$ that 
are mapped to each other. Since the hyper-K\"ahler metric shall be $\rho$-invariant,
we have to blow up the singularities that are described by the $\widetilde{w}'_{2i+1}$, 
too, if we blow up the singularities that are described by the $\widetilde{w}'_{2i}$. 
With help of this geometric interpretation, we are able to formulate our corollary.

\begin{Co} 
\label{K3SingCor}
Let $(r,a,\delta)\in\mathbb{N}\times\mathbb{N}_0\times \{0,1\}$ be a 
triple such that there exists a K3 surface with a simple non-symplectic involution 
with invariants $(r,a,\delta)$. Moreover, let $S$ be the K3 surface 
from Theorem \ref{K3SingThm} that has $3$ points with $A_1$-singularities 
and $2$ points with $E_8$-singularities and let $\rho$ be the
non-symplectic involution from the same theorem. 
Moreover, let $G_1,\ldots,G_{k_1}$ be the connected components of 
$3A_1\cup 2E_8$ that are mapped to itself by $\rho$ and let $G'_1,\ldots,
G'_{k_2}$ be a set of connected components that are not invariant under
$\rho$ such that  

\[
G_1\cup\ldots\cup G_{k_1} \cup G'_1 \cup \ldots \cup G'_{k_2}
 \cup \rho(G'_1) \cup \ldots \cup \rho(G'_{k_2}) = 3A_1\cup 2E_8
\]

Finally, let $\widetilde{G}_1,\ldots,\widetilde{G}_{l_1}$ be  
connected Dynkin diagrams that can be obtained by deleting some nodes 
of $G_1\cup\ldots\cup G_{k_1}$ and let  $\widetilde{G}'_1,\ldots,
\widetilde{G}'_{l_2}$  be connected Dynkin diagrams that can 
be obtained by deleting some nodes of $G'_1\cup\ldots\cup G'_{k_2}$. 
Then there exists a K3 surface with a 
hyper-K\"ahler metric that admits an isometric non-symplectic involution
with invariants $(r,a,\delta)$ that has $l_1$ singular points of type 
 $\widetilde{G}_1,\ldots,\widetilde{G}_{l_1}$ and $2l_2$ singular points
of type $\widetilde{G}_1,\ldots,\widetilde{G}_{l_2}$.
\end{Co}

\begin{Ex}
Let $(r,a,\delta)=(10,10,0)$. We recall that this is the case where the 
fixed locus is empty. It is possible to choose the marking such that 
$\rho$ acts as $\rho_1^7\oplus \rho_2^4$ on $L$. More explicitly, 
we have

\begin{equation*}
\rho(u_j^i) =u_j^{3-i}\:,\quad
\rho(u_j^3) = -u_j^3\:,\quad
\rho(v_k^i) =v_k^{3-i}
\end{equation*} 

for all $i,j\in\{1,2\}$ and $k\in\{1,\ldots,8\}$. $\rho$ interchanges the two 
Dynkin diagrams of type $E_8$ and two of the Dynkin diagrams of type
$A_1$. The third Dynkin diagram of type $A_1$ is preserved by $\rho$ 
since $\rho(\widetilde{w}_3) = \rho(u^3_1 - u^3_2) = -\widetilde{w}_3$ 
and $-\widetilde{w}_3$ is another root of the lattice $A_1$. We delete the 
node from $E_8$ that is connected to three other nodes. The remaining 
diagram is of type $A_1\cup A_2\cup A_4$. Corollary \ref{K3SingCor}
guarantees that there exists a singular K3 surface with a non-symplectic 
involution $\rho$ with invariants $(10,10,0)$ that has $5$ singular points 
of type $A_1$, $2$ of type $A_2$ and $2$ of type $A_4$. Both points of 
type $A_2$ and $A_4$ are mapped by $\rho$ to each other. Moreover, 
there exist $2$ points with $A_1$-singularities that are mapped to $2$ 
other points with $A_1$-singularities and one point $p\in S$ with an 
$A_1$-singularity is fixed by $\rho$. 

We remove the singularity at $p$ such that $\rho:L\rightarrow L$ is still 
induced by a non-symplectic involution. This is only possible if we add a 
term $\lambda\widetilde{w}_3$ to $x$ or $y$. Afterwards, $\widetilde{w}_3$ 
is not contained in the Picard lattice anymore and thus it does not correspond 
to a complex curve on the K3 surface. Technically speaking, the family of K3 
surfaces with $x_t :=x + t\widetilde{w}_3$ defines a one-parameter family of 
hyper-K\"ahler metrics that converges to the singular one, but our construction 
is not a resolution in the sense of algebraic geometry. Since $\rho$ acts 
orientation-reversing on the 2-sphere that represents $\widetilde{w}_3$, our 
example does not contradict the fact that an involution with invariants $(10,10,0)$ 
of a smooth K3 surface does not have any fixed points. 
\end{Ex}

\section{K3 surfaces with two involutions}
\label{TwoInvolutions}

Finally, we study K3 surfaces with two commuting involutions that are
non-symplectic with respect to two anti-commuting complex structures. 
As before, let $x,y,z\in L_{\mathbb{R}}$ be the images of the $3$
K\"ahler classes with respect to the marking. We denote the two non-symplectic
involutions by $\rho^1,\rho^2:S\rightarrow S$. Without loss of generality, $\rho^1$
and $\rho^2$ shall act on the K\"ahler classes as

\begin{equation}
\label{KahlerRelations}
\begin{array}{ccc}
\rho^1(x)=-x & \rho^1(y)=-y & \rho^1(z)=z \\
\rho^2(x)=-x & \rho^2(y)=y & \rho^2(z)=-z \\
\end{array}
\end{equation}

The composition $\rho^1\rho^2$ is a third involution that 
satisfies 

\begin{equation}
\rho^1\rho^2(x)=x \quad \rho^1\rho^2(y)=-y \quad \rho^1\rho^2(z)=-z
\end{equation}

There is a straightforward method to construct pairs $(\rho^1,\rho^2)$
with the above properties. Let $(r_i,a_i,\delta_i)$ with $i=1,2$ be triples
of invariants that belong to non-symplectic involutions such that the direct sum 
$L^1\oplus L^2$ of the fixed lattices can be primitively embedded into
$L$. Theorem \ref{LatticeEmbThm} guarantees that this is possible if

\[
r_1+r_2\leq 11 \quad\text{or}\quad r_1+r_2+a_1+a_2 < 22
\]

Since the embedding is primitive, there exists a basis

\begin{equation}
\label{L1L2basis}
(u_1,\ldots,u_{r_1},v_1,\ldots,v_{r_2},w_1,\ldots,w_{22-r_1-r_2})
\end{equation}

of $L$ such that $(u_1,\ldots,u_{r_1})$ is a basis of $L^1$ and 
$(v_1,\ldots,v_{r_2})$ is a basis of $L^2$. The maps $\rho^1$ 
and $\rho^2$ that are defined by  

\[
\begin{array}{ccc}
\rho^1(u_i)= u_i & \rho^1(v_i)=-v_i & \rho^1(w_i)= -w_i \\
\rho^2(u_i)= -u_i & \rho^2(v_i)= v_i & \rho^2(w_i)= -w_i \\
\end{array}
\]

commute and they induce non-symplectic involutions with respect to
suitable complex structures. Since $L^1$ and $L^2$ are hyperbolic 
lattices and $L$ has signature $(3,19)$, the lattice $(L^1\oplus L^2)^{\perp}
= \text{span}_{\mathbb{Z}}(w_1,\ldots,w_{22-r_1-r_2})$ is hyperbolic, too.
Therefore, it is possible to choose $z\in L^1_{\mathbb{R}}$, 
$y\in L^2_{\mathbb{R}}$ and $x\in (L^1_{\mathbb{R}}\oplus 
L^2_{\mathbb{R}})^{\perp}$ such that $x^2=y^2=z^2>0$. Since
the three lattices are pairwise orthogonal, we have $x\cdot y=
x\cdot z = y\cdot z = 0$ automatically. Moreover, $x$, $y$ and $z$
satisfy the relations (\ref{KahlerRelations}). All in all, we have
constructed a K3 surface with a hyper-K\"ahler structure and 
two commuting involutions that are non-symplectic with respect
to different complex structures. Since the sets of all positive elements
in $L^1_{\mathbb{R}}$, $L^2_{\mathbb{R}}$ or $(L^1_{\mathbb{R}}
\oplus L^2_{\mathbb{R}})^{\perp}$ are open, we can choose

\[
x= \sum_{i=1}^{22 - r_1 - r_2} \gamma_i w_i \quad
y = \sum_{i=1}^{r_2} \beta_i v_i \quad
z = \sum_{i=1}^{r_1} \alpha_i u_i  
\]

such that 

\[
(\alpha_1,\ldots,\alpha_{r_1},\beta_1,\ldots,\beta_{r_2},
\gamma_1,\ldots,\gamma_{22-r_1-r_2})
\]

is $\mathbb{Q}$-linearly independent. Since any $d\in L$ has integer
coefficients with respect to the basis (\ref{L1L2basis}), this condition 
guarantees that there exists no $d\in L$ with $d^2=-2$ and $x\cdot d= 
y\cdot d = z\cdot d = 0$. Therefore, it is possible to choose $S$ as a 
smooth K3 surface. All in all, we have proven the following theorem. 

\begin{Th}
Let $(r_1,a_1,\delta_1),(r_2,a_2,\delta_2) \in \mathbb{N}\times \mathbb{N}_0
\times \{0,1\}$ such that there exist non-symplectic involutions with invariants
$(r_i,a_i,\delta_i)$. Moreover, we assume that $r_1+r_2\leq 11$ or $r_1+r_2+a_1+a_2 
< 22$. In this situation, there exists a smooth K3 surface $S$ with a hyper-K\"ahler
structure that admits two commuting involutions $\rho^1$ and $\rho^2$ that are 
non-symplectic with respect to different complex structures $I_1$ and $I_2$ with 
$I_1I_2=-I_2I_1$ and have invariants $(r_1,a_1,\delta_1)$ and $(r_2,a_2,\delta_2)$. 
\end{Th}

\begin{Rem}
An important step in Kovalev's and Lee's construction of $G_2$-manifolds 
\cite{KoLe} is to find two K3 surfaces $S_1$ and $S_2$ with non-symplectic
involutions $\rho^1$ and $\rho^2$ and a so called matching. A matching 
is defined as an isometry $f:S_1\rightarrow S_2$ such that  

\[
f^{\ast} \omega_{I_2} = \omega_{J_1}\:,
\qquad
f^{\ast} \omega_{J_2} = \omega_{I_1}\:,
\qquad 
f^{\ast} \omega_{K_2} = -\omega_{K_1}\:,
\]

where $\omega_{I_k}$, $\omega_{J_k}$ and $\omega_{K_k}$ are
the three K\"ahler forms on $S_k$. Let $S$ be a K3 surface with two
involutions that satisfy (\ref{KahlerRelations}). If we choose the triple of 
complex structures on $S$ first as $(I,J,K)$ and then as $(J,I,-K)$, the 
identity map becomes a matching. Therefore, the above theorem 
shows that a matching exists if the invariants of $\rho^1$ and $\rho^2$
satisfy $r_1+r_2\leq 11$ or $r_1+r_2+a_1+a_2 < 22$. This fact is also
shown in \cite{KoLe}. 
\end{Rem}

We are especially interested in constructing K3 surfaces with ADE-sin\-gularities 
that admit a pair of commuting involutions with the above properties. Unfortunately, 
it is hard to tell how the set $D$ that describes the singular set can look like in this
general situation. The reason for this is that we have an existence result for 
the basis (\ref{L1L2basis}) but no further information. If $\rho^1$ and
$\rho^2$ are simple, it is possible to choose the hyper-K\"ahler structure
such that we can determine $D$ explicitly. 

Therefore, we assume from now on that $\rho^1$ and $\rho^2$ are
simple. Let $\phi_k:H^2(S,\mathbb{Z})\rightarrow L$ with $k=1,2$ be
markings such that $\rho^k(w_i)=\pm w_j$. The matrix representation 
of the pull-back map $\rho^1:H^2(S,\mathbb{Z}) \rightarrow 
H^2(S,\mathbb{Z})$ with respect to the basis $\phi_1^{-1}(w_i)$ 
is a convenient matrix whose columns are unit vectors multiplied with 
$\pm 1$. If $\phi_2\neq\phi_1$, the matrix representation of $\rho^2$
with respect to $\phi_1^{-1}(w_i)$ may be more complicated. Since that case
is rather difficult to handle, we restrict ourselves to the case $\phi_1=\phi_2$.

Up to conjugation, $\rho^k:L\rightarrow L$ can be written 
as $\rho^{i_k}_1 \oplus \rho^{j_k}_2$, where $\rho^{i_k}_1:3H 
\rightarrow 3H$ and $\rho^{j_k}_2:2(-E_8)\rightarrow 2(-E_8)$
are two of the maps that we have defined in Section \ref{NonsympSection}.
By a direct calculation we see that $\rho^{j_1}_2$ and $\rho^{j_2}_2$
commute if and only if  $(j_1,j_2)\notin \{(2,4),(4,2)\}$. By adjusting the 
marking $\phi_1$, we can assume that the restriction of $\rho^1$ to
$2(-E_8)$ actually is one of the maps $\rho^{j_1}_2$ with 
$j_1\in\{1,\ldots,4\}$. Nevertheless, the restriction of $\rho^2$ 
may be a conjugate of a map $\rho^{j_2}_2$ such that we still have 
$\rho^2(w_i)=\pm w_j$ for $i\in\{7,\ldots,22\}$. As we have remarked 
in Section \ref{NonsympSection}, the only additional possibilities 
for $\rho^2|_{2(-E_8)}$ are 

\[
\rho^2|_{2(-E_8)}(x_1,x_2)=(x_1,-x_2)
\]

if $j_2=2$ or 

\[
\rho^2|_{2(-E_8)}(x_1,x_2)=(-x_2,-x_1)
\]

if $j_2=4$. If we take account of these additional possibilities, it is 
still not possible that $\rho^1|_{2(-E_8)}$ and $\rho^2|_{2(-E_8)}$
commute if $(j_1,j_2)\in \{(2,4),(4,2)\}$. Nevertheless, this idea will 
be helpful in the next case. Let $i_1,i_2\in \{1,\ldots,7\}$. First, we 
assume that $i_1,i_2\neq 7$. We see that $\rho^{i_1}_1$ and 
$\rho^{i_2}_1$ always commute, since the smaller matrix blocks

\[
\left(\,
\begin{array}{cc}
1 & 0 \\
0 & 1 \\
\end{array}
\,\right)\:, \quad
\left(\,
\begin{array}{cc}
-1 & 0 \\
0 & -1 \\
\end{array}
\,\right)\:, \quad
\left(\,
\begin{array}{cc}
0 & 1 \\
1 & 0 \\
\end{array}
\,\right)\:, \quad
\left(\,
\begin{array}{cc}
0 & -1 \\
-1 & 0 \\
\end{array}
\,\right) 
\]

commute pairwisely. In Section \ref{OneInvolution} we have defined a hyper-K\"ahler
structure by 

\[
x:= u_1^2+u_2^2\:,\quad
y:= u_1^3+u_2^3\:,\quad
z:= u_1^1+u_2^1\:.
\]

The involution $\rho^1$ preserves $z$ and acts as $-1$ on $x$ and $y$. Unfortunately,
the same is true for $\rho^2$, although $\rho^2$ should preserve $y$ and act as $-1$ 
on $x$ and $z$. In order to solve this problem, we conjugate $\rho^{i_2}_1$ by the 
map $\tau:3H\rightarrow 3H$ that is defined by

\[
\tau(u_k^l) := u_k^{4-l} \quad\forall k\in\{1,2\},l\in\{1,2,3\}\:.
\]

In other words, we permute the first and the third block of the matrices
that define $\rho_1^{i_2}$. We obtain a map that is still an isometry of $3H$
and maps any $w_i$ to a $\pm w_j$.  After this conjugation, $\rho^{i_1}_1$
and $\rho^{i_2}_1$ still commute and the maps $\rho^1,\rho^2:L \rightarrow
L$ satisfy the relations (\ref{KahlerRelations}).  

If $i_1=7$ and $i_2\in\{1,\ldots,6\}$, $\rho_1^{i_1}$ has a $4\times 4$-block 
in the upper left corner that interchanges $H_1$ and $H_2$. Therefore, 
$\rho^{i_1}_1$ and $\tau^{-1}\rho^{i_2}_1\tau$ commute if and only if the 
last two $2\times 2$-blocks of $\rho^{i_2}_1$ are the same. This is the case 
for all values of $i_2$ except $2$ and $5$.  We consider the second 
hyper-K\"ahler structure from Section \ref{OneInvolution} that is defined by

 \[
x:= u_1^1+u_2^1 - u_1^2 - u_2^2\:,\quad
y:= \sqrt{2}(u_1^3+u_2^3)\:,\quad
z:= u_1^1+u_2^1 + u_1^2 + u_2^2 \:.
\]

After a short calculation, we see that $\rho^1$ and $\rho^2$ satisfy
the relations (\ref{KahlerRelations}) again. All in all, we 
have proven the following sufficient condition for the existence of a pair 
$(\rho^1,\rho^2)$ of non-symplectic involutions. 

\begin{Th}
\label{Z22actionThm}
Let $(i_1,j_1),(i_2,j_2)\in\{1,\ldots,7\}\times\{1,\ldots,4\}$ such that $(j_1,j_2)\notin
\{(2,4),(4,2)\}$ and $(i_1,i_2)\notin \{(2,7),(5,7),(7,2),(7,5),(7,7)\}$. Moreover, let \linebreak
$(r_k,a_k,\delta_k)$ with $k\in\{1,2\}$ be the triples of invariants that characterise the 
non-symplectic involutions that act as $\rho_1^{i_k} \oplus \rho_2^{j_k}$
on $L$. In this situation, there exists a possibly singular K3 surface $S$ that 
admits two commuting involutions $\rho^1$ and $\rho^2$ that are non-symplectic 
with respect to different complex structures $I_1$ and $I_2$ with $I_1I_2=-I_2I_1$ 
and have invariants $(r_1,a_1,\delta_1)$ and $(r_2,a_2,\delta_2)$. 
\end{Th}

\begin{Rem}
The above theorem yields $320$ different sets $\{(r_k,a_k,\delta_k) | k\in\{1,2\} \}$
of invariants of pairs $(\rho^1,\rho^2)$ with the desired properties. We remark that 
our result is mainly an existence theorem. For one set of invariants there may exist more 
than one pair of simple non-symplectic involutions with the same invariants. If we 
choose for example $\rho_2^{j_2}$ as one of the maps (\ref{rho2alt}) or modify 
$\rho_1^{i_2}$ by permuting the three summands $H_1$, $H_2$ and $H_3$, 
we could easily obtain further examples with the same invariants but a different 
action of $\mathbb{Z}^2_2$ on $L$. Since we have restricted ourselves to the case 
that the $\rho^k$ are simple and both markings $\phi_k:H^2(S,\mathbb{Z})
\rightarrow L$ are the same, it is even possible that examples with further sets 
of invariants exist. The investigation of these questions is beyond the scope of
this paper. 
\end{Rem}

Since the hyper-K\"ahler structure on $S$ that we have introduced in the proof of
the theorem is the same as in Section \ref{OneInvolution}, we immediately obtain 
the following corollary.

\begin{Co}
\label{Z22actionCor}
In the situation of the above theorem, $S$ can be chosen as a K3 surface
that has $3$ singular points with $A_1$-singularities and $2$ singular points
with $E_8$-singularities.
\end{Co}

Our next step is to investigate if there exist K3 surfaces with further kinds of 
singularities that admit involutions $\rho^1$ and $\rho^2$ with the same 
properties as in Theorem \ref{Z22actionThm}. Let $(\widetilde{w}_i)_{i=1,\ldots,19}$ 
be the basis (\ref{PicardMaxBasis}) of the lattice that we have introduced in 
(\ref{PicardMax}). We recall that $\widetilde{w}_i^2 = -2$ for all $i$ and that the 
$\widetilde{w}_i$ correspond to the nodes of the Dynkin diagram $3A_1\cup 2E_8$. 
The involutions $\rho^1$ and $\rho^2$ generate a group that is isomorphic to 
$\mathbb{Z}^2_2$. We denote the span of the orbit of $\widetilde{w}_i$ by $W_i$. 
The dimension of $W_i$ is either $1$, $2$ or $4$. For the same 
reasons as in Section \ref{OneInvolution}, $\mathbb{Z}^2_2$ acts on
$3A_1\cup 2E_8$ and maps connected components to connected 
components. Since $3A_1\cup 2E_8$ does not contain 
$4$ components of the same type, the dimension of $W_i$ has to be $1$ or $2$. 
We call a $\widetilde{w}_i$ of type 

\begin{itemize}
    \item $(1,1)$ if $\rho^1(\widetilde{w}_i)=\rho^2(\widetilde{w}_i)=\widetilde{w}_i$,
    \item $(1,-1)$ if $\rho^1(\widetilde{w}_i)=\widetilde{w}_i$ and $\rho^2(\widetilde{w}_i)
    \neq \widetilde{w}_i$, 
    \item $(-1,1)$ if $\rho^1(\widetilde{w}_i)\neq \widetilde{w}_i$ and $\rho^2(\widetilde{w}_i)= 
    \widetilde{w}_i$,
    \item $(-1,-1)$ if $\rho^1(\widetilde{w}_i)\neq \widetilde{w}_i$ and $\rho^2(\widetilde{w}_i)
    \neq \widetilde{w}_i$.
\end{itemize}

Since $\rho^1$ and $\rho^2$ are involutions that preserve exactly one positive vector, 
their eigenvalues are precisely $1$ and $-1$. Moreover, they commute and therefore 
we have a decomposition

\[
L_{\mathbb{R}} = V_{1,1} \oplus V_{1,-1} \oplus V_{-1,1} \oplus V_{-1,-1} 
\]

where 

\[
V_{\epsilon_1,\epsilon_2} = \{ v\in L_{\mathbb{R}} | \rho^1(v)=\epsilon_1 v, 
 \rho^2(v)=\epsilon_2 v\}\:.
\]

We have $x\in V_{-1,-1}$, $y\in V_{-1,1}$ and $z\in V_{1,-1}$. If $\widetilde{w}_i$ is of type
$(1,-1)$, we define a $\widetilde{w}'_i\in L_{\mathbb{R}}$ by

\[
\widetilde{w}'_i = 
\begin{cases}
\widetilde{w}_i & \text{if}\:, \rho^2(\widetilde{w}_i)=-\widetilde{w}_i \\
\widetilde{w}_i - \widetilde{w}_j & \text{if}\:, \rho^2(\widetilde{w}_i)=\widetilde{w}_j\:\text{with}\:i\neq j \\
\widetilde{w}_i + \widetilde{w}_j & \text{if}\:, \rho^2(\widetilde{w}_i)=-\widetilde{w}_j\:\text{with}\:i\neq j \\
\end{cases}
\]  

If $\widetilde{w}_i$ is of type $(-1,1)$, we define $\widetilde{w}'_i$ analogously but replace $\rho^2$ by 
$\rho^1$. Finally, if $\widetilde{w}_i$ is of type $(-1,-1)$, we define

\[
\widetilde{w}'_i = 
\begin{cases}
\widetilde{w}_i & \text{if}\:, \rho^1(\widetilde{w}_i)=\rho^2(\widetilde{w}_i)=-\widetilde{w}_i \\
\widetilde{w}_i - \widetilde{w}_j & \text{if}\:, \rho^1(\widetilde{w}_i)=-\widetilde{w}_i\:\text{and}
\:\rho^2(\widetilde{w}_i)=\widetilde{w}_j\:\text{with}\:i\neq j \\
\widetilde{w}_i + \widetilde{w}_j & \text{if}\:, \rho^1(\widetilde{w}_i)=-\widetilde{w}_i\:\text{and}
\:\rho^2(\widetilde{w}_i)=-\widetilde{w}_j\:\text{with}\:i\neq j \\
\widetilde{w}_i - \widetilde{w}_j & \text{if}\:, \rho^2(\widetilde{w}_i)=-\widetilde{w}_i\:\text{and}
\:\rho^1(\widetilde{w}_i)=\widetilde{w}_j\:\text{with}\:i\neq j \\
\widetilde{w}_i + \widetilde{w}_j & \text{if}\:, \rho^2(\widetilde{w}_i)=-\widetilde{w}_i\:\text{and}
\:\rho^1(\widetilde{w}_i)=-\widetilde{w}_j\:\text{with}\:i\neq j \\
\end{cases}
\] 

Since $\dim{W_i}\neq 4$, these are the only possibilities that can happen for 
a $\widetilde{w}_i$ of type $(-1,-1)$. By our construction $\widetilde{w}'_i\in 
V_{\epsilon_1,\epsilon_2}$ if $\widetilde{w}'_i$ is of type $(\epsilon_1,\epsilon_2)$. 
We choose arbitrary subsets 

\[
\begin{array}{rcl}
P & \subseteq & \{1\leq i\leq 19 | \text{$\widetilde{w}_i$ is of type $(-1,-1)$}\} \\
Q & \subseteq & \{1\leq i\leq 19 | \text{$\widetilde{w}_i$ is of type $(-1,1)$}\} \\
R & \subseteq & \{1\leq i\leq 19 | \text{$\widetilde{w}_i$ is of type $(1,-1)$}\} \\
\end{array}
\]

such that for any pair $(i,j)$ with $i\neq j$ from one the three sets we have
$W_i\cap W_j= \{0\}$. Let $(x,y,z)$ be the triple of K\"ahler 
classes that determines the hyper-K\"ahler structure with $3$ $A_1$- and $2$ 
$E_8$-singularities. We define a new hyper-K\"ahler structure by 

\[
\begin{array}{rcl}
x' & = & \mu x + \sum_{i\in P} \alpha_i \widetilde{w}'_i \\
y' & = & \nu y + \sum_{i\in Q} \beta_i \widetilde{w}'_i \\
z' & = & \lambda z + \sum_{i\in R} \gamma_i \widetilde{w}'_i \\
\end{array}
\]  

The coefficients in the above definition are chosen such that 
\begin{enumerate}
    \item the family that consists of $1$, the $\alpha_i$, the $\beta_i$ and the
    $\gamma_i$ is $\mathbb{Q}$-linearly independent,
    \item $x'^2=y'^2=z'^2 > 0$. 
\end{enumerate}

The hyper-K\"ahler structure that is defined by $x'$, $y'$ and $z'$ still satisfies 
the equation (\ref{KahlerRelations}). Moreover, the set $D$ that determines
the number and type of the singular points can be obtained from $3A_1\cup 2E_8$
by deleting all nodes that correspond to an element of the $\mathbb{Z}^2_2$-orbit 
of an $i\in P\cup Q\cup R$. In other words, we have constructed a (partial) resolution 
of the singularities that is still invariant under $\mathbb{Z}^2_2$. We remark that 
in general there is a minimal singularity that cannot be resolved without destroying 
the $\mathbb{Z}^2_2$-symmetry. Its Dynkin diagram is given by all $i$ such that
$\widetilde{w}_i$ is invariant under $\mathbb{Z}^2_2$. If we add a multiple of such an 
$\widetilde{w}_i$ to $x$, $y$ or $z$, we obtain a new hyper-K\"ahler structure that 
no longer satisfies (\ref{KahlerRelations}). All in all, we have proven the following 
theorem. 

\begin{Th} 
\label{Z22actionThm2}
Let $S$ be one of the K3 surfaces from Theorem \ref{Z22actionThm}
that 

\begin{enumerate}
    \item admits a pair $(\rho^1,\rho^2)$ of commuting simple involutions 
    that are non-symplectic with respect to two complex structures $I_1$ and
    $I_2$ with $I_1I_2=-I_1I_2$ and
    
    \item has $3$ points with $A_1$-singularities and $2$ points 
    with $E_8$-singularities.
\end{enumerate}    

$\rho^1$ and $\rho^2$ generate a group that is isomorphic to
$\mathbb{Z}^2_2$ and acts on the Dynkin diagram
$3A_1\cup 2E_8$. Let $M$ be a $\mathbb{Z}^2_2$-invariant subset 
of the nodes of $3A_1\cup 2E_8$ such that no node from $M$
corresponds to a $\widetilde{w}_i\in L$ that is fixed by
$\mathbb{Z}^2_2$. In this situation, there exists a K3 surface $S'$
that  

\begin{enumerate}
    \item admits a pair of commuting simple involutions that are 
    non-symplectic with respect to two complex structures $I'_1$ and
    $I'_2$ with $I'_1I'_2=-I'_1I'_2$ and whose invariants $(r_i,a_i,\delta_i)$
    are the same as of $\rho^i$ and
    
    \item whose singular set is described by the Dynkin diagram that we obtain
    by deleting the set $M$ of nodes from $3A_1\cup 2E_8$.  
\end{enumerate}  

In particular, $S'$ can be chosen as a smooth K3 surface if there is no 
$\widetilde{w}_i$ that is fixed by $\mathbb{Z}^2_2$
\end{Th} 

\begin{Ex}
Let $\rho^i:L\rightarrow L$ with $i=1,2$ be the lattice isometries that act as 
the identity on $H_i\oplus 2(-E_8)$ and as $-1$ on the other two summands
that are isometric to $H$. $\rho^1$ and $\rho^2$ commute and are both of 
type $\rho^1_1 \oplus \rho^1_2$. Corollary \ref{Z22actionCor} guarantees
that there exists a K3 surface with two $E_8$- and three $A_1$-singularities 
and two non-symplectic involutions that correspond to $\rho^1$ and $\rho^2$.    
Theorem \ref{Z22actionThm2} allows us to resolve one or two of the 
$A_1$-singularities, but the two $E_8$-singularities and the third of the 
$A_1$-singularities cannot be resolved without destroying the invariance of
the hyper-K\"ahler metric with respect to $\rho^1$ and $\rho^2$.   
\end{Ex}


\begin{thebibliography}{999}
    \bibitem{Ach} Acharya, B.S.: M-Theory, $G_2$-manifolds and Four Dimensional Physics.
    Classical and Quantum Gravity, Vol.19 (2002) No. 22., 5619 - 5653.  

    \bibitem{AchGuk} Acharya, B.S.; Gukov, S.: M-theory and singularities of exceptional 
    holonomy manifolds. Phys. Rep. 392 (2004), No. 3, 121 - 189. 

    \bibitem{Anderson1} Anderson, M.T.: Moduli spaces of Einstein metrics on 4-manifolds.
    Bulletin of the AMS 21 (1989), No. 2, 275 - 279. 
    
    \bibitem{Anderson2} Anderson, M.T.: The $L^2$-structure of moduli spaces of Einstein
    metrics on $4$-manifolds. Geometric and Functional Analysis 2 (1992), No. 1, 29 - 89. 

    \bibitem{ArSaTa} Artebani, M.; Sarti, A.; Taki, S.: K3 surfaces with non-symplectic 
    automorphisms of prime order. With an appendix by Shigeyuki Kond$\bar{\text{o}}$. 
    Math. Z. 268 (2011), No. 1-2, 507 - 533.

    \bibitem{BHPV}  Barth, W.; Hulek, K.; Peters, C.; van de Ven, A.: Compact complex
    surfaces. Second Enlarged Edition. Springer-Verlag, Berlin Heidelberg, 2004.

    \bibitem{Do} Dolgachev, I.: Integral quadratic forms: Applications to algebraic geometry
    (after V. Nikulin). Bourbaki Seminar Vol. 1982/83. Ast\'{e}risque 105 - 106 (1983), 
    251 - 278. 

    \bibitem{Do2} Dolgachev, I. V.: Mirror symmetry for lattice polarized K3 surfaces. Algebraic 
    geometry. 4. J. Math. Sci. 81 (1996), No. 3, 2599 - 2630.

    \bibitem{Ebeling} Ebeling, W.: Lattices and codes. A course partially based on lectures by 
    F. Hirzebruch. Second revised edition. Advanced Lectures in Mathematics. Friedr. Vieweg \& Sohn,
    Braunschweig, 2002.

    \bibitem{Joyce} Joyce, D.: Compact manifolds with special holonomy. Oxford Mathematical 
    Monographs. Oxford University Press, Oxford, 2000.

    \bibitem{KarJoy} Karigiannis, S.; Joyce, D.: A new construction of compact torsion-free 
    $G_2$-manifolds by gluing families of Eguchi-Hanson spaces. Preprint. arXiv:1707.09325 
    [math.DG]

    \bibitem{Koba} Kobayashi, R.; Todorov, A.N.: Polarized period map for generalized K3
    surfaces and the moduli of Einstein metrics. Tohoku Math. Journ. 39 (1987),
    341 - 363. 

    \bibitem{KoLe} Kovalev, A.; Lee, N.-H.: K3 surfaces with non-symplectic involution and 
    compact irreducible $G_2$-manifolds. Math. Proc. Cambridge Philos. Soc. 151 (2011), 
    No. 2, 193 - 218.

    \bibitem{Nikulin1} Nikulin, V.V.: Finite groups of automorphisms of K\"ahlerian K3
    surfaces. Trans. Moscow Math. Soc. 2 (1980), 71 - 135. 

    \bibitem{Nikulin2} Nikulin, V.V.: Integer symmetric bilinear forms and some of their 
    applications. Math. USSR Izvestia 14 (1980), 103 - 167. 

    \bibitem{Nikulin3} Nikulin, V. V.: Factor groups of groups of automorphisms of hyperbolic 
    forms by subgroups generated by 2-reflections. Algebro-geometric applications.  
    J. Soviet Math. 22 (1983), 1401- 1476.
    
    \bibitem{ReidHabil} Reidegeld, Frank: $G_2$-orbifolds with ADE-singularities.
    Habilitation thesis. Fakult\"at f\"ur Mathematik, TU Dortmund, 2017. Online available:
    http://dx.doi.org/10.17877/DE290R-18940
\end{thebibliography}
\end{document}